%%%%%%%%%%%%%%%%   Geometry and Topology: 2003-4.tex  %%%%%%%%%%%%%
%%%%        
%%%%    The smooth Whitehead spectrum of a point at odd regular primes
%%%%             
%%%%                            John Rognes 
%%%%  
%%%%                Published in Volume 7(2003) pages 155-184
%%%%
%%%%                   Publication date 17 February 2003
%%%%
%%%%                        This is a plain TeX file
%%%%
%%%%
%%%%%%%%%%%%%%%%%%                                   %%%%%%%%%%%%%%%%%%%
%%%%%%%%%%%%%%%%%%%%%%%%%%%%%%%%%%%%%%%%%%%%%%%%%%%%%%%%%%%%%%%
%%%%%%%%%%%             gtmacros.tex            %%%%%%%%%%%%%%%
%%%%%%%%%%%             version 1.6             %%%%%%%%%%%%%%% 
%
%                       Colin Rourke   
%
%
%    These macros are recommended for use by authors submitting articles   
%    to Geometry and Topology or to Algebraic and Geometric Topology.  
%    They are intended to be used with plain TeX. Each macro is described 
%    briefly to make it clear how to use it (or to modify it to achieve
%    different results).  If you modify this file then please change its
%    name.  If you modify this file and use the modified file to 
%    format an article for submission to Geometry and Topology or
%    Algebraic and Geometric Topology, then please paste the modified
%    file into your main TeX file.  Do not submit it as a separate file.
%      
%    Instructions on using these macros are also given in  gtmacins.tex  
%    or  gtmacins.ps  or .pdf  available on the gt www pages or by 
%    anonymous ftp from the gt/info/macros directory.
%
%
\magnification=\magstephalf      % Sets default point size to 11pt.
%
%  Basic layout parameters :
%
\vsize=7.5truein                 % Sets text height to 7.5 inches.
\hsize=5.2truein                 % Sets text width to 5.2 inches.
\newskip\stdskip                 % standard vertical space
\stdskip=6pt plus3pt minus3pt    % (slightly more stretchy
\medskipamount=\stdskip          % than the usual \medskip)
\parindent=0pt                   % Paragraphs are non-indented with
\parskip=\stdskip                % a little space between paragraphs. 
\abovedisplayskip=\stdskip       %  Reduces the space
\belowdisplayskip=\stdskip       %  around displays.
\mathsurround=0.75pt             % Gives a little extra space around maths.
\overfullrule=0pt                %  Prevents black boxes
%
%   The following macro is for principal paragraph breaks ie
%   a paragraph break with a slightly larger space :
%
\def\ppar{\par\goodbreak\vskip 8pt plus 4pt minus 4pt}     
%
%  The standard horizontal space for theorems, labels etc :
%
\def\stdspace{\hskip 0.75em plus 0.15em\ignorespaces}
\let\qua\stdspace % useful abbreviation (3/4 of a quad)
%
%%%%%%%%%%%%%%            FONT MACROS            %%%%%%%%%%%%
%
%           The following font macros define the AMS symbol 
%           and Euler-Fraktal fonts for use in text and
%           mathematics with appropriate size changes.
%           They also define two new control sequences  
%           \small  and  \large  (similar to those built
%           into LaTeX) which change the size of all fonts 
%           both in text and maths.  \small  is 10% smaller 
%           than normal and  \large  30% bigger.  The strange
%           size of the \large text fonts (10pt scaled 1315)
%           is because these macros are intended to be used
%           at \magstephalf.  The result is 10pt scaled 1440
%           (\magstep2) which is a standard font size.  If
%           you are borrowing these macros to use them at
%           another basic  \magnification, then you will
%           probably need to change 1315 to 1200 in the eleven
%           places marked ** below.  \large  will then be
%           20% bigger than normal.  Note that at \magstephalf
%           all the fonts come out roughly one point larger
%           than their size as defined in these macros.
%
%           The size-changing macros are based on Knuth's
%           \ninepoint and \eightpoint macros.
%
%
%    The macros are laid out in a way which makes it clear how to
%    add futher fonts (or delete unavailable ones) and how to add
%    further size changes.
%
%    First comes a definition of  \hexnumber  which is needed to
%    refer to font families whose family number is not known :
%
\def\hexnumber#1{\ifcase#1 0\or 1\or 2\or 3\or 4\or 5\or 6\or 7\or 8\or
 9\or A\or B\or C\or D\or E\or F\fi}
%
%     Next we define the AMS symbol-a fonts at 13,10,9,7,6,5 points
%
\font\thirtnmsa=msam10 scaled 1315    %%% **  see note above 
\font\tenmsa=msam10          \font\ninemsa=msam9
\font\sevenmsa=msam7         \font\sixmsa=msam6
\font\fivemsa=msam5
%%%%%%  (add further sizes here if you need them)
%
%    and the standard size family for these fonts
%
\newfam\msafam                  \textfont\msafam=\tenmsa
\scriptfont\msafam=\sevenmsa    \scriptscriptfont\msafam=\fivemsa
\edef\hexa{\hexnumber\msafam}        %  The msa family is  \fam\hexa
\def\msa{\fam\msafam\tenmsa}         %  \msa  switches to this family
%
%    Repeat these steps for the AMS symbol-b fonts
%
\font\thirtnmsb=msbm10 scaled 1315   %%%  ** see note above
\font\tenmsb=msbm10      \font\ninemsb=msbm9
\font\sevenmsb=msbm7     \font\sixmsb=msbm6
\font\fivemsb=msbm5
%%%%%%  (add further sizes here if you need them)
%
\newfam\msbfam                   \textfont\msbfam=\tenmsb       
\scriptfont\msbfam=\sevenmsb     \scriptscriptfont\msbfam=\fivemsb
\edef\hexb{\hexnumber\msbfam}    %  The msb family is \fam\hexb  
\def\msb{\fam\msbfam\tenmsb}     %  \msb switches to this family
%
%        Repeat for the Euler-Fraktal fonts 
%
\font\thirtneufm=eufm10 scaled 1315   %%% **  see note above 
\font\teneufm=eufm10                 \font\nineeufm=eufm9
\font\seveneufm=eufm7                \font\sixeufm=eufm6
\font\fiveeufm=eufm5
%%%%%%  (add further sizes here if you need them)
%
\newfam\eufmfam                    \textfont\eufmfam=\teneufm
\scriptfont\eufmfam=\seveneufm     \scriptscriptfont\eufmfam=\fiveeufm
\edef\hexf{\hexnumber\eufmfam}      % The Euler-Fraktal family is
\def\frak{\fam\eufmfam\teneufm}     % \fam\hexf and \frak switches to this
%
%%%  Add further fonts families here (using the same format) if you need
%    them.  The def of hexnumber is optional (it is only used for
%    \mathchardef 's).
%
%      Now we need to define the standard fonts (which are
%      already defined at 10,7 and 5 point) at 13,9 and 6 point:
%
%      Roman fonts:
\font\thirtnrm=cmr10 scaled 1315    %%%  ** see note above
\font\ninerm=cmr9                   \font\sixrm=cmr6   
%%%%%%  (add further sizes here if you need them)
%
%      Math italic fonts
\font\thirtni=cmmi10 scaled 1315    %%%  ** see note above 
\font\ninei=cmmi9                   \font\sixi=cmmi6  
%%%%%%  (add further sizes here if you need them)
%
%     Symbol fonts
\font\thirtnsy=cmsy10 scaled 1315   %%%  ** see note above
\font\ninesy=cmsy9                  \font\sixsy=cmsy6  
%%%%%%  (add further sizes here if you need them)
%
%     Bold face
\font\thirtnbf=cmbx10 scaled 1315   %%%  ** see note above 
\font\ninebf=cmbx9                  \font\sixbf=cmbx6  
%%%%%%  (add further sizes here if you need them)
%
%     The maths extension font (only defined at text size)
%
\font\thirtnex=cmex10 scaled 1315   %%%  ** see note above
\font\nineex=cmex9                  
%%%%%%  (add further sizes here if you need them)
%
%     Finally three fonts (text italic, slanted and typewriter type)
%     which are also only defined at text size
%
\font\thirtnit=cmti10 scaled 1315  %%%  ** see note above 
\font\nineit=cmti9                  
%%%%%%  (add further sizes here if you need them)
%
\font\thirtnsl=cmsl10 scaled 1315  %%%  ** see note above 
\font\ninesl=cmsl9                  
%%%%%%  (add further sizes here if you need them)
%
\font\thirtntt=cmtt10 scaled 1315  %%%  ** see note above 
\font\ninett=cmtt9                  
%%%%%%  (add further sizes here if you need them)
%
%
%     Now come the two main macros.  What  \small  does is to
%     change all the families of fonts from normal size which is
%     10,7,5  (ie 10pt text, 7pt subscript, 5pt subsubscript)
%     to 9,6,5.  \large  similarly changes to  13,9,7.  To make
%     other size changing macros, choose your three sizes, add
%     font size definitions if necessary and make the obvious changes
%     to one of these macros.  Change  \normalbaselineskip  and
%     \strutbox  dimensions to appropriate sizes as well.  To
%     add further fonts, insert them in each macro, using the
%     AMS fonts as a model.
%      
%
\def\small{%
%
%   redefine the sizes of the roman fonts :
%
\textfont0=\ninerm \scriptfont0=\sixrm \scriptscriptfont0=\fiverm
\def\rm{\fam0\ninerm}%       % ( \rm  sets \ninerm  in text mode
%                            %  and \fam0 in math mode)
%
%   and the math italic fonts :
%
\textfont1=\ninei \scriptfont1=\sixi \scriptscriptfont1=\fivei
%
%   and the symbol fonts :
%
\textfont2=\ninesy \scriptfont2=\sixsy \scriptscriptfont2=\fivesy
%
%   There is only one math extension font :
%
\textfont3=\nineex \scriptfont3=\nineex \scriptscriptfont3=\nineex
%
%   Next the bold font (named rather than numbered) :
%
\textfont\bffam=\ninebf \scriptfont\bffam=\sixbf
\scriptscriptfont\bffam=\fivebf \def\bf{\fam\bffam\ninebf}%
%
%   and the three text-only fonts : 
%
\textfont\itfam=\nineit \def\it{\fam\itfam\nineit}%
\textfont\slfam=\ninesl \def\sl{\fam\slfam\ninesl}%
\textfont\ttfam=\ninett \def\tt{\fam\ttfam\ninett}%
%
%   Now the three new families of AMS fonts :
%
%   AMS symbol-a
%
\textfont\msafam=\ninemsa \scriptfont\msafam=\sixmsa
\scriptscriptfont\msafam=\fivemsa \def\msa{\fam\msafam\ninemsa}%         
%
%   AMS symbol-b
%
\textfont\msbfam=\ninemsb \scriptfont\msbfam=\sixmsb
\scriptscriptfont\msbfam=\fivemsb \def\msb{\fam\msbfam\ninemsb}%         
%
%   Euler-Fraktal font
%
\textfont\eufmfam=\nineeufm  \scriptfont\eufmfam=\sixeufm
\scriptscriptfont\eufmfam=\fiveeufm \def\frak{\fam\eufmfam\nineeufm}%
%
%%%  Add further fonts families here if you need them.
%
%    Reset \normalbaselineskip and \strubox :
%
\normalbaselineskip=11pt%
\setbox\strutbox=\hbox{\vrule height8pt depth3pt width0pt}%
%
%    Set \normalbaselines and \rm (roman) as defaults :
%
\normalbaselines\rm
%
%    Reset some of the basic vertical skips:
%
\stdskip=4pt plus2pt minus2pt    
\medskipamount=\stdskip          
\parskip=\stdskip                
\abovedisplayskip=\stdskip       
\belowdisplayskip=\stdskip       
\def\ppar{\par\goodbreak\vskip 6pt plus 3pt minus 3pt}%     
%
%   And finally reset the size of section heads (see below):
%
\def\section##1{\global\advance\sectionnumber by 1
\vskip-\lastskip\penalty-800\vskip 20pt plus10pt minus5pt 
\egroup{\bf\number\sectionnumber\quad##1}\bgroup\small         
\vskip 6pt plus3pt minus3pt
\nobreak\resultnumber=1}%      % Reset resultnumber at start of section
}    %%%   End of  \small  macro      
%
%   Two useful abbreviations to keep track of \small material:
\def\beginsmall{\bgroup\small}
\let\endsmall\egroup
%
%
%    The \large  macro is similar (comments abbreviated):
%
%
\def\large{%
\textfont0=\thirtnrm \scriptfont0=\ninerm \scriptscriptfont0=\sevenrm
\def\rm{\fam0\thirtnrm}%
\textfont1=\thirtni \scriptfont1=\ninei \scriptscriptfont1=\seveni
\textfont2=\thirtnsy \scriptfont2=\ninesy \scriptscriptfont2=\sevensy
\textfont3=\thirtnex \scriptfont3=\thirtnex \scriptscriptfont3=\thirtnex
\textfont\bffam=\thirtnbf \scriptfont\bffam=\ninebf
\scriptscriptfont\bffam=\sevenbf \def\bf{\fam\bffam\thirtnbf}%
\textfont\itfam=\thirtnit \def\it{\fam\itfam\thirtnit}%
\textfont\slfam=\thirtnsl \def\sl{\fam\slfam\thirtnsl}%
\textfont\ttfam=\thirtntt \def\tt{\fam\ttfam\thirtntt}%
%   AMS symbol-a  :
\textfont\msafam=\thirtnmsa \scriptfont\msafam=\ninemsa
\scriptscriptfont\msafam=\sevenmsa \def\msa{\fam\msafam\thirtnmsa}%         
%   AMS symbol-b  :
\textfont\msbfam=\thirtnmsb \scriptfont\msbfam=\ninemsb
\scriptscriptfont\msbfam=\sevenmsb \def\msb{\fam\msbfam\thirtnmsb}%         
%   Euler-Fraktal font :
\textfont\eufmfam=\thirtneufm  \scriptfont\eufmfam=\nineeufm
\scriptscriptfont\eufmfam=\seveneufm \def\frak{\fam\eufmfam\teneufm}%
%%%% Add further fonts families here if you need them.
%   Reset \normalbaselineskip and \strubox and initialise :
\normalbaselineskip=16pt%
\setbox\strutbox=\hbox{\vrule height11.5pt depth4.5pt width0pt}%
\normalbaselines\rm}%
\let\Large\large   %  for compatibility with latex
%
%   The next two lines define commonly used switches for
%   blackboard bold (\Bbb) and gothic type (\goth).  The   
%   \Bbb  switch is set to work in the same way as in amstex
%   and switches only the next character to blackboard bold.
\def\Bbb#1{{\msb#1}}

%
%   To use the new AMS fonts you can either use the control
%   sequences \msa, \msb (alias \Bbb) and \frak (alias \goth) eg :
\def\N{{\frak N}}

   % see the msam font table
%
%   or, more generally, make \mathchardef's (cf Knuth p155) eg :
\mathchardef\plussquare="0\hexa01
\mathchardef\nge="3\hexb0B
\mathchardef\maltesecross="0\hexa7A
\mathchardef\del="0\hexf01
%
%   or you can use the amstex names for all the new symbols by
%   inserting the line  \input amsnames  in your file directly
%   after \input gtmacros. 
%   This presupposes that you have collected a copy of the file
%   amsnames.tex  from the  gt/info/macros  ftp directory.
%
%
%   Finally we need a small capital font (for author(s)) :
%
\font\sc=cmcsc10
%
%%%%%%%%%%%%%%%%%       END OF FONT MACROS     %%%%%%%%%%%%%
%
%
%                 Knuth's \square macro :
%
\def\sqr#1#2{{\vcenter{\vbox{\hrule  height.#2truept
	\hbox{\vrule width.#2truept height#1truept 
	\kern#1truept \vrule width.#2truept}
	\hrule height.#2truept}}}}
\def\sq{\sqr55}    %   A small square for end-of-proofs. 
%                  %   (Define other size squares by varing the
%                  %   the two numbers.)
%
%
%      Style macros for section heads, theorem statements etc :
%   
%
\newcount\sectionnumber            %%%  Allocate registers to take
\newcount\resultnumber             %%%  section and result numbers.
\sectionnumber=0\resultnumber=1    %%%  Set these registers to 0 and 1
%
%   The \section macro produces a \large bold faced section heading
%   numbered to the left.  Pagebreaks are encouraged before the
%   start of the section and discouraged directly after the heading.
%   Typical use  \section{First steps}  with typical result :
%
%    1  First Steps     (set bold and \large)
%
\def\section#1{\global\advance\sectionnumber by 1
\xdef\nextkey{\number\sectionnumber}%      (used by the \key macro)
\vskip-\lastskip\penalty-800\vskip 20pt plus10pt minus5pt 
{\large\bf\number\sectionnumber\quad#1}         
\vskip 8pt plus4pt minus4pt
\nobreak\resultnumber=1}      % Reset resultnumber at start of section
%
%
%
%   Next a macro to set subheadings (like the  \section  macro
%   but without the number, with less space and set in standard size).
%
%   Typical use :  \sh{Example formats}
%
\def\sh#1{\vskip-\lastskip\ppar{\bf #1}\par\nobreak\medskip}         
%
%   The \proc ... \endproc macros ("proclaim") are for setting theorems, 
%   lemmas, conjectures etc with automatic numbering.  Typical use :    
%  
%    \proc{Theorem}Every lemon is yellow.\endproc
%
%   Typical result :
%     
%    Theorem 3.4  Every lemon is yellow.   

%   (with Theorem 3.4 set bold and a \stdspace of space before the 
%   statement set in slanted type).
%
\def\proc#1{\xdef\nextkey{\number\sectionnumber.\number\resultnumber}%
\vskip-\lastskip\ppar\bf%
\noindent#1\ \number\sectionnumber.\number\resultnumber
\stdspace\sl\global\advance\resultnumber by 1\ignorespaces}
 
%
%  The \prf ... \endprf macros are for setting proofs.  The code
%  for \prf includes the code for \endproc, so there is no need to
%  type \endproc if the theorem is followed immediatedly by a proof.
%
                            %  For start of proofs  
\def\qed{\hfill$\sq$\par\goodbreak\rm}   %  For end (or absence) of proofs
                 %  extra vertical space)  
        %  For start of proof with alternative name
              %  \endproof is an alias for \endprf
%
%   Typical uses :    
%  
%    \proc{Theorem}Every lemon is yellow. \qed\endproc
%
%    \proc{Theorem}Every lemon is yellow.
%    \prf Use your eyes. \endprf
%
%    \proc{Theorem}Every lemon is yellow.
%    \proof{Proof of theorem} Use your eyes. \endprf
%
%   The next macro is a variant of the \proc macro.  It has
%   exactly the same result except that it omits the number.
%
%   Typical use :  
%    
%    \proclaim{Conjecture}Some oranges are yellow.\endproc
%
\def\proclaim#1{\vskip-\lastskip\ppar\bf%
\noindent#1\stdspace\sl\ignorespaces} 
\let\endproclaim\endproc
%
%   The next macro is a further variant for remarks, definitions etc.   
%   It omits the number and does not switch on slanted type.  
%  
%   Typical use :
%
%    \rk{Remark}Some lemons are thick-skinned.\endrk
%
\def\rk#1{\vskip-\lastskip\ppar{\bf #1}\stdspace\ignorespaces}                

%
%   The next macro is for numbering equations etc, \label  produces the 
%   correct number  x.y  and advances the resultnumber register
%
%   Typical use :
%
%     $$fx=7\eqno{\bf\label}$$
%
%   result :
%
%                           fx = 7                           3.5
%
\def\label{\xdef\nextkey{\number\sectionnumber.\number\resultnumber}%
\number\sectionnumber.\number\resultnumber
\global\advance\resultnumber by 1}
%
%
%
%   The next macros are to automate external references.  To use them 
%   type \reflist ..... \endreflist near the beginning of your paper, 
%   where  .... is the list of references in alphabetical order 
%   and in  the form  \key{KEY}  reference    where "KEY" is a 
%   string of characters which reminds you of the reference.   
%   Separate  references with a blank line or a \par.   Eg 
%
%     \reflist
%
%     ..... more references ....
%
%     \key{Kn-84} {\bf D Knuth}, {\it The TeXbook}, Addison--Wesley (1984)
%
%     ..... more references ....
%
%     \endreflist
%
%   Then type  \references  where you wish the references to be printed
%   (normally near the end of the paper).  To refer to Knuth type
%   for example    see Knuth [\ref{Kn-84}, page 320]   and the correct
%   numerical reference will be printed.  Edit the \references macro
%   to change the formatting (if desired).
%   There is an alternative \refkey for \key, provided your KEY contains
%   only letters.  The syntax is:
%
%     \reflist
%
%     ..... more references ....
%
%     \refkey\Knuth  {\bf D Knuth}, {\it The TeXbook}, Addison--Wesley (1984)
%
%     ..... more references ....
%
%     \endreflist
%
%   \key{Knuth}  has exactly the same maening as \refkey\Knuth and you
%   can mix the two syntaxes if you want.  But \refkey\Kn-84
%   would not work.  It would set Kn as the KEY and -84 would get printed!
%
\newcount\refnumber              %  Register for reference numbers
\refnumber=1                     %  set initially to 1.
\long\def\reflist#1\endreflist{%
\long\def\thereflist{#1}{\def\refkey##1##2\par{\xdef##1{\number\refnumber}%
\global\advance\refnumber by 1}%
\def\key##1##2\par{\expandafter\xdef%
\csname##1\endcsname{\number\refnumber}%
\global\advance\refnumber by 1}#1\par}}
\long\def\references{%
\penalty-800\vskip-\lastskip\vskip 15pt plus10pt minus5pt 
{\large\bf References}\ppar %`References' is set \large bold with space around.
{\leftskip=25pt\frenchspacing    % The list of references is set 
\small\parskip=3pt plus2pt       % \small  with small spaces between,
\def\refkey##1##2\par{\noindent  % numbers in [,]'s and set just to the
\llap{[##1]\stdspace}\ignorespaces##2\par}         % left of a 25pt margin.
\def\key##1##2\par{\noindent  
\llap{[\ref{##1}]\stdspace}\ignorespaces##2\par}  
\def\,{\thinspace}\thereflist\par}}
%
%   Next a footnote macro (with automatic numbering) which sets the
%   footnote  \small.
%
%   Typical use :
%         ..... are yellow.\fnote{By yellow here we mean Britsh
%    Standard colour BS3320.} 
%
\newcount\footnotenumber         % Register for footnote number
\footnotenumber=1                % set initially to 1
\def\fnote#1{\xdef\nextkey{\number\footnotenumber}%
{\small\ifnum\footnotenumber>9\parindent=14pt%
\else\parindent=10pt\fi\footnote{$^{\number\footnotenumber}$}%
{\hglue-5pt#1}\global\advance\footnotenumber by 1}}
%
%
%   Next macros for handling figures with automatic numbering (using 
%   TeX's \midinsert to float the figure to a suitable place).
%   
%   The \figure ... \endfigure macro centres the figure and adds
%   an automatically numbered label  Figure XX  after it.
%
%   If you have a caption, then type \caption{caption text} 
%   somewhere between \figure and \endfigure.  The macro
%   will then add  Figure XX: caption text  after the figure.
%
%   If you want an unnumbered or uncentred figure, then use TeX's raw 
%       \midinsert Figure instructions \endinsert  
%   and if you want a numbered figure label in the same style then
%   use \caption{caption text} outside of  \figure ... \endfigure.
%
%   If you need just the label Figure XX  outside of  \figure ... \endfigure
%   then type  \figurelabel .
%
\newcount\figurenumber          % register for figure number
\figurenumber=1                 % set initially to 1
\def\caption#1{\xdef\nextkey{\number\figurenumber}%
\cl{\small Figure \number\figurenumber: #1}%
\global\advance\figurenumber by 1}
\def\figurelabel{\xdef\nextkey{\number\figurenumber}%
\cl{\small Figure \number\figurenumber}%
\global\advance\figurenumber by 1}
\long\def\figure#1\endfigure{{\xdef\nextkey{\number\figurenumber}%
\let\captiontext\relax\def\caption##1{\xdef\captiontext{##1}}%
\midinsert\cl{\ignorespaces#1\unskip\unskip\unskip\unskip}\vglue6pt\cl{\small 
Figure \number\figurenumber\ifx\captiontext\relax\else: \captiontext
\fi}\endinsert\global\advance\figurenumber by 1}}
%
%   Macros for self-correcting internal references.
%
%   There are two macros  \key{KEY}  and  \ref{KEY} .
%
%   The \key macro sets up KEY as a key for whatever number is 
%   being referenced and the \ref macro converts the KEY into 
%   that number.  Type \key after a  \section or \proc or 
%   \label or \fnote or \figure or \caption or \figurelabel .
%
%   Example:
%
%       \section{Introduction}\key{intro}
%       \proc{Theorem}\key{MainTh}Lemons are yelloy\endproc
%       Here we follow\fnote{Follow in the sense of Dickens}
%       \key{Dickens-note}the crowd ....  
%
%       In section \ref{intro}
%       we stated theorem \ref{mainTh} and noted (see footnote 
%       \ref{Dickens-note}) ...
%
\def\nextkey{??}   %  initialise \nextkey (which is reset by all the
%                     numbering macros)
%
\def\key#1{\expandafter\xdef\csname #1\endcsname{\nextkey}}
\def\ref#1{\expandafter\ifx\csname #1\endcsname\relax
\immediate\write16{Reference {#1} undefined}??\else
\csname #1\endcsname\fi}
%
%   Note:  If the KEY contains only letters then \KEY has exactly the
%   same meaning as \ref{KEY} so in the example you could have:
%
%       In section \intro\ we ....
%
%   The \key will work at any time after the macro which sets the
%   number, provided no other macro which sets a number has been used. 
%
%   Macros for forward references:
%              =======
%   The \key \ref macros ONLY work for backwards references.  If you  
%   want to use forwards references, then type \useforwardrefs  near
%   the beginning of your file.  The KEY's are then stored in an
%   auxiliary  .ref  file and you then suffer the same disadvantage as
%   when using LaTeX that you must TeX the file twice to get
%   the references correct.
%
%   To use a forward ref type \ref{KEY}.  (You can type the
%   alternative  \KEY  but you'll get an error on first TeX'ing 
%   if the \KEY is not yet defined.) 
%
%   The macro also allows external references to be listed at the end 
%   of the file (if you wish to).  (Indeed they can be typed anywhere
%   before the \references command.)  You can combine the reference list
%   and the \references command by typing the references (using the
%   same syntax as before) between the commands \biblio and \endbiblio 
%   (don't type \references or they'll be printed twice).
%
\newread\gtinfile
\newwrite\gtreffile
\def\useforwardrefs{
\openin\gtinfile\jobname.ref
\ifeof\gtinfile
\closein\gtinfile
\immediate\write16{No file \jobname.ref}
\else
\closein\gtinfile
\input \jobname.ref
\fi
\immediate\openout\gtreffile \jobname.ref
%
%   Adapt \key :
%
\def\key##1{{\def\\{\noexpand}%
\expandafter\xdef\csname ##1\endcsname{\nextkey}%
\immediate\write\gtreffile{\\\expandafter\\\def\\\csname ##1\\\endcsname%
{\nextkey}}}}
%
%  Adapt macros for external references:  
%
\long\def\reflist##1\endreflist{%
\long\def\thereflist{##1}{\def\refkey####1####2\par{\xdef####1{%
\number\refnumber}{\def\\{\noexpand}\immediate\write\gtreffile
{\\\def\\####1{\number\refnumber}}}\global\advance\refnumber by 1}%
\def\key####1####2\par{\expandafter\xdef%
\csname####1\endcsname{\number\refnumber}%
{\def\\{\noexpand}\immediate\write\gtreffile
{\\\expandafter\\\def\\\csname ####1\\\endcsname{\number\refnumber}}}
\global\advance\refnumber by 1}##1\par}}
\long\def\biblio##1\endbiblio{\reflist##1\endreflist\references}%
%
%  Adapt obselete key macros (\numkey, \seckey and \figkey):
%
\def\numkey##1{{\def\\{\noexpand}%
\xdef##1{\number\sectionnumber.\number\resultnumber}
\immediate\write\gtreffile{\\\def\\##1%
{\number\sectionnumber.\number\resultnumber}}}}
\def\seckey##1{{\def\\{\noexpand}\xdef##1{\number\sectionnumber}
\immediate\write\gtreffile{\\\def\\##1{\number\sectionnumber}}}}
\def\figkey##1{\xdef##1{\number\figurenumber}%
{\def\\{\noexpand}\immediate\write\gtreffile%
{\\\def\\##1{\number\figurenumber}}}
\number\figurenumber\global\advance\figurenumber by 1}
}   %  end of \useforwardrefs
%
%
%   The next five macros are obselete and have been superseeded by
%   the general \key macro above.  They are included merely to 
%   maintain backward compatibility for the package:
%
%
\def\figkey#1{\xdef#1{\number\figurenumber}%
\number\figurenumber\global\advance\figurenumber by 1}
\def\fig#1#2\endfig{%
\midinsert\cl{#2}\vglue6pt\cl{\small Figure #1}\endinsert}
\def\newfig{\number\figurenumber\global\advance\figurenumber by 1}
\def\numkey#1{\xdef#1{\number\sectionnumber.\number\resultnumber}}
\def\seckey#1{\xdef#1{\number\sectionnumber}}
%
%   End of obselete macros.
%
%
%   The next macro is a version of the verbatim macro given by Knuth.
%
%   This macro produces a "verbatim" printout of
%   any ASCII string which does not contain the symbol "
%   (TeX files do not usually contain " 's).
%   More precisely, everything between consecutive pairs
%   of " 's is printed verbatim in the typewriter font cmtt.
%   For an explanation of how the macro works, see Knuth pp 420-1.
%
%   There are two switches: \verb (which switches the macro on)
%   and \brev which switches the macro off (the default).  When
%   the macro is switched off the symbol " has its usual 
%   meaning for TeX.  To use the macro, type \verb before use
%   and the use " to switch verbatim on and off.  Be careful
%   not to use " for any other purpose.  There is no need to
%   switch the macro off again unless you need to use " for
%   some other purpose (eg making  \mathchardef 's).  Note 
%   that the macro MUST BE OFF before inputting  amsnames.tex .
%
%   Whether the macro is on or off you can always use the
%   control sequence \dq (double quote) for " e.g.
%   \mathchardef\sum=\dq1350  is perfectly valid.
%   The control sequence \ttq is an abbreviation for
%   {\tt\dq}.  Thus "\ttq" will produce " (in cmtt)
%   inside a verbatim quote.
%
%
   %  define a code for " so it can be used when \verb is on
  %  code for " in cmtt
%
\def\verb{\catcode`\"=\active}       %  The main
\def\brev{\catcode`\"=12}            %  switches.
\brev                                %  Prime switches and
\verb                                %  switch on.
{\obeyspaces\gdef {\ }}              
{\catcode`\`=\active\gdef`{\relax\lq}}
\def"{%
\begingroup\baselineskip=12pt\def\par{\leavevmode\endgraf}%
\tt\obeylines\obeyspaces\parskip=0pt\parindent=0pt%
\catcode`\$=12\catcode`\&=12\catcode`\^=12\catcode`\#=12%
\catcode`\_=12\catcode`\~=12%
\catcode`\{=12\catcode`\}=12\catcode`\%=12\catcode`\\=12%
\catcode`\`=\active\let"\endgroup}
\brev      %   Finally switch the macro off (for safety)
%
%   Macros for itemised lists.   Typical use :
%    
%    \items
%    \item{(i)}Colours must be defined.
%    \item{(ii)}Colour cards may not be cited.
%    \enditems
%
%   Result :
%
%    (i)  Colours must be defined. 
%   (ii)  Colour cards may not be cited.
%
%
           % Start of itemised list         
         % end of itemised list   
\def\item#1{\par\leavevmode\llap{#1\stdspace}%
\ignorespaces}                             % labelled item
               % bulleted item.
%
%   The \quote ... \endquote macros are for typesetting quotations :
%

%
%   A few useful abbreviations :
%
    %  Colon with correct spacing for maps.
\def\np{\vfil\eject}         %  Forced page break (new page).
\def\nl{\hfil\break}         %  New line.
\def\cl{\centerline}         %  Centerline
\def\gt{{\mathsurround=0pt\it $\cal G\mskip-2mu$eometry \&\ 
$\cal T\!\!$opology}}        %  The journal title in recommended style
    %  for monographs
\def\agt{{\mathsurround=0pt\it$\cal A\mskip-.7mu$lgebraic \&\ 
$\cal G\mskip-2mu$eometric $\cal T\!\!$opology}}  % AGT
%
%    Finally some macros for automatic title page or header generation.
%    To use them type your header information using the following  
%    example as a guide :
%
%    Note that \\ is used as standard separator (for lines in \title and
%    \address, between authors and between email addresses or URL's)
%    and that \email, \url and \secondaddress are optional.
%

% Example:  \title{A short spoof paper\\with a two-line title}
% =======   \authors{Albert Einstein\\Leonardo da Vinci}
%           \address{IAS\\Princeton}\secondaddress{Renaissance\\Venice}
%           \email{ae@ias.princeton.edu\\ldv@ren.ven.hist}
%           \abstract 
%           A short spoof paper with a very short abstract.
%           \endabstract 
%           \primaryclass{00-01, 00-02}\secondaryclass{68-00, 68-01}
%           \keywords{Short, spoof, paper}
%           \maketitlepage
%
%
%    The title page or header will then be generated automatically.
%
%
%    Define the various ingredients of the title page:
%
\def\title#1{\def\thetitle{#1}}
\def\shorttitle#1{\def\theshorttitle{#1}}
\def\author#1{\edef\previousauthors{\theauthors}
 \ifx\theauthors\relax\def\theauthors{#1}\else
 \def\theauthors{\previousauthors\par#1}\fi}

        % aliases
%
\def\address#1{\edef\previousaddresses{\theaddress}
 \ifx\theaddress\relax\def\theaddress{#1}\else
 \def\theaddress{\previousaddresses\par\vskip 2pt\par#1}\fi}
                             % alias
\def\secondaddress#1{\edef\previousaddresses{\theaddress}
 \ifx\theaddress\relax\def\theaddress{#1}\else
 \def\theaddress{\previousaddresses\par{\rm and}\par#1}\fi}   

\def\email#1{\edef\previousemails{\theemail}
 \ifx\theemail\relax\def\theemail{#1}\else
 \def\theemail{\previousemails\hskip 0.75em\relax#1}\fi}
  % aliases
\def\secondemail#1{\edef\previousemails{\theemail}
 \ifx\theemail\relax\def\theemail{#1}\else
 \def\theemail{\previousemails\hskip 0.75em{\rm and}\hskip 0.75em
 \relax#1}\fi}
\def\url#1{\edef\previousurls{\theurl}
 \ifx\theurl\relax\def\theurl{#1}\else
 \def\theurl{\previousurls\hskip 0.75em\relax#1}\fi}
      % aliases
\def\secondurl#1{\edef\previousurls{\theurl}
 \ifx\theurl\relax\def\theurl{#1}\else
 \def\theurl{\previousurls\hskip 0.75em{\rm and}\hskip 0.75em
 \relax#1}\fi}
\long\def\abstract#1\endabstract{\long\def\theabstract{#1}}
\def\primaryclass#1{\def\theprimaryclass{#1}}
\let\subjclass\primaryclass                        % alias
\def\secondaryclass#1{\def\thesecondaryclass{#1}}
\def\keywords#1{\def\thekeywords{#1}}
%
%  Set \\ to \par and title page items to \relax to initialise macros :
%
\let\\\par\let\thetitle\relax\let\theshorttitle\relax
\let\theauthors\relax\let\theshortauthors\relax
\let\theaddress\relax\let\theshortaddress\relax
\let\theemail\relax\let\theurl\relax
\let\theabstract\relax\let\theprimaryclass\relax
\let\thesecondaryclass\relax\let\thekeywords\relax
%
%
%
%   Basic title page layout (edit this macro if you
%   wish to adjust the title page layout) :
%
\long\def\maketitlepage{    % start of definition of \maketitlepage

\vglue 0.2truein   % top margin

% title :
%
{\parskip=0pt\leftskip 0pt plus 1fil\def\\{\par\smallskip}{\large
\bf\thetitle}\par\medskip}   

\vglue 0.15truein 

% authors :
%
{\parskip=0pt\leftskip 0pt plus 1fil\def\\{\par}{\sc\theauthors}
\par\medskip}%
 
\vglue 0.1truein 

% address(es) email's and URL's (with switches to detect whether the
% optional items have been used) :
%
{\small\parskip=0pt
{\leftskip 0pt plus 1fil\def\\{\par}{\sl\theaddress}\par}
\ifx\theemail\relax\else  % email address?
\vglue 5pt \def\\{\stdspace{\rm and}\stdspace} 
\cl{Email:\stdspace\tt\theemail}\fi
\ifx\theurl\relax\else    % URL given?
\vglue 5pt \def\\{\stdspace{\rm and}\stdspace} 
\cl{URL:\stdspace\tt\theurl}\fi\par}

\vglue 7pt 

{\bf Abstract}

\vglue 5pt

\theabstract

\vglue 7pt 

{\bf AMS Classification numbers}\quad Primary:\quad \theprimaryclass\par

Secondary:\quad \thesecondaryclass

\vglue 5pt 

{\bf Keywords:}\quad \thekeywords

\np  % page break at the end of the title page

}    % end of definition of \maketitlepage
%
%    % \makeshorttitle (for general preprints) doesn't take a new page
%
\long\def\makeshorttitle{    % start of definition of \makeshorttitle

%\vglue 0.2truein   % top margin

% title :
%
{\parskip=0pt\leftskip 0pt plus 1fil\def\\{\par\smallskip}{\large
\bf\thetitle}\par\medskip}   

\vglue 0.05truein 

% authors :
%
{\parskip=0pt\leftskip 0pt plus 1fil\def\\{\par}{\sc\theauthors}
\par\medskip}%
 
\vglue 0.03truein 

% address(es) email's and URL's (with switches to detect whether the
% optional items have been used) :
%
{\small\parskip=0pt
{\leftskip 0pt plus 1fil\def\\{\par}{\sl\ifx\theshortaddress\relax
\theaddress\else\theshortaddress\fi}\par}
\ifx\theemail\relax\else  % email address?
\vglue 5pt \def\\{\stdspace{\rm and}\stdspace} 
\cl{Email:\stdspace\tt\theemail}\fi
\ifx\theurl\relax\else    % URL given?
\vglue 5pt \def\\{\stdspace{\rm and}\stdspace} 
\cl{URL:\stdspace\tt\theurl}\fi\par}

\vglue 10pt 

% abstract and classification numbers (with switches):

{\small\leftskip 25pt\rightskip 25pt{\bf Abstract}\stdspace\theabstract

{\bf AMS Classification}\stdspace\theprimaryclass
\ifx\thesecondaryclass\relax\else; \thesecondaryclass\fi\par
{\bf Keywords}\stdspace \thekeywords\par}
\vglue 7pt
}    % end of definition of \makeshorttitle
\let\maketitle\makeshorttitle        %% alias
%
%    %%%% \makeagttitle (for AGT) similar to \makeshorttitle but
%         with addresses omitted (they go at the end)
%
%%%% publication info and test defaults:

\def\volumenumber#1{\def\thevolumenumber{#1}}
\def\volumeyear#1{\def\thevolumeyear{#1}}
\def\pagenumbers#1#2{\def\startpage{#1}\def\finishpage{#2}}
\def\published#1{\def\publishdate{#1}}
\def\received#1{\def\receiveddate{#1}}
\def\revised#1{\def\reviseddate{#1}}
\let\reviseddate\relax
%% Defaults for authors to use to check layout
\volumenumber{X}
\volumeyear{20XX}
\pagenumbers{1}{XXX}
\published{XX Xxxember 20XX}

\long\def\makeagttitle{   %%% start of definition of \makeagttitle
\agt\hfill      %   Journal title (top left) 
%   logo placeholder (top right)
\hbox to 60truept{\vbox to 0pt{\vglue -14truept{\bf [Logo here]}\vss}\hss}
\break
{\small Volume \thevolumenumber\ (\thevolumeyear)
\startpage--\finishpage\nl
Published: \publishdate}

\vglue .2truein

% title
{\parskip=0pt\leftskip 0pt plus 1fil\def\\{\par\smallskip}{\large
\bf\thetitle}\par\medskip}   
\vglue 0.05truein 

% authors :
%
{\parskip=0pt\leftskip 0pt plus 1fil\def\\{\par}{\sc\theauthors}
\par\medskip}%
 
\vglue 0.03truein 

%  abstract and classification numbers:

{\small\leftskip 25truept\rightskip 25truept{\bf Abstract}\stdspace\theabstract

{\bf AMS Classification}\stdspace\theprimaryclass
\ifx\thesecondaryclass\relax\else; \thesecondaryclass\fi\par
{\bf Keywords}\stdspace \thekeywords\par}\vglue 7truept

}   %%%% end of definition of \makeagttitle

%%%%% Macro to typeset addresses (typically at the end of the paper)

\def\Addresses{\bigskip
{\small \parskip 0pt \leftskip 0pt \rightskip 0pt plus 1fil \def\\{\par}
\sl\theaddress\par\medskip \rm Email:\stdspace\tt\theemail\par
\ifx\theurl\relax\else\smallskip \rm URL:\stdspace\tt\theurl\par\fi}}

\def\agtart{%   Full mock-up of AGT article style (for authors to test with)
%  get print centerpage:
\hoffset 14truemm
\voffset 31truemm
\font\phead=cmsl9 scaled 950
\font\pnum=cmbx10 scaled 913
\font\pfoot=cmsl9 scaled 950
%  headline and footline
\headline{\vbox to 0pt{\vskip -4.5mm\line{\small\phead\ifnum
\count0=\startpage ISSN numbers are printed here
\hfill {\pnum\folio}\else\ifodd\count0\def\\{ }% 
\ifx\theshorttitle\relax\thetitle\else\theshorttitle\fi\hfill{\pnum\folio}
\else\def\\{ and }{\pnum\folio}\hfill\ifx\theshortauthors\relax\theauthors
\else\theshortauthors\fi\fi\fi}\vss}}
\footline{\vbox to 0pt{\vglue 0mm\line{\small\pfoot\ifnum\count0=\startpage
Copyright declaration is printed here\hfill\else
\agt, Volume \thevolumenumber\ (\thevolumeyear)\hfill\fi}\vss}}
%  force \agttitle
\let\maketitle\makeagttitle\let\makeshorttitle\makeagttitle}

%%%
%%%  This version of  gtoutput.tex  is intended to finish formatting
%%%  papers published in Geometry & Topology and stored in the
%%%  arXiv.   All versions of  gtoutput.tex  are copyright 
%%%  GT Publications and are to be used _only_ for formatting
%%%  the officially published version of G&T papers.
%%%
%%%
%%%                                             Colin Rourke  14.9.2000
%%%
%%%  To create header file  head.xxx  comment out the first \endinput

%  test for latex or plain tex
\def\ifplaintex{\expandafter\ifx\csname documentclass\endcsname\relax}

%  get print centerpage:

\ifplaintex 
\hoffset 14truemm
\voffset 31truemm
\else
\headsep 23pt
\footskip 35pt
\hoffset -4truemm
\voffset 12.5truemm
\fi

%  load pictex if not already loaded :
\expandafter\ifx\csname beginpicture\endcsname\relax
\expandafter\ifx\csname documentclass\endcsname\relax
\input pictex \else\font\fiverm=cmr5
\input prepictex \input pictex \input postpictex \fi\fi

\def\gt{{\mathsurround=0pt\it $\cal G\mskip-2mu$eometry \&\ 
$\cal T\!\!$opology}}        %  journal title in recommended style

\def\gtp{{\mathsurround=0pt\it $\cal G\mskip-2mu$eometry \&\ 
$\cal T\!\!$opology $\cal P\!$ublications}}  % GT publications

%  define the various new ingredients of the title page 

\def\lognumber#1{\def\thelognumber{#1}}
\def\volumenumber#1{\def\thevolumenumber{#1}}
\def\papernumber#1{\def\thepapernumber{#1}}
\def\volumeyear#1{\def\thevolumeyear{#1}}

\def\pagenumbers#1#2{\def\startpage{#1}\def\finishpage{#2}}
\def\published#1{\def\publishdate{#1}}
\def\proposed#1{\def\theproposer{#1}}
\def\seconded#1{\def\theseconders{#1}}
\def\received#1{\def\receiveddate{#1}}
\def\revised#1{\def\reviseddate{#1}}
\def\accepted#1{\def\accepteddate{#1}}

\long\def\asciiabstract#1{\long\def\theasciiabstract{#1}}
\def\asciikeywords#1{\def\theasciikeywords{#1}}

\def\shorttitle#1{\def\theshorttitle{#1}}

%  initialise

\let\\\par\let\thelognumber\relax
\let\thevolumenumber\relax\let\thepapernumber\relax
\let\thevolumeyear\relax\let\thesamplenumber\relax\let\startpage\relax
\let\finishpage\relax\let\publishdate\relax\let\receiveddate\relax
\let\reviseddate\relax\let\accepteddate\relax\let\theasciititle\relax
\let\theasciiauthors\relax
\let\theasciiabstract\relax\let\theasciikeywords\relax
\let\theasciiemail\relax\let\theshortauthors\relax\let\theshorttitle\relax

\long\def\maketitlep{   % start of definition of \maketitlep

\count0=\startpage

\gt\hfill      %   Journal title (top left) 
%    Logo (top right) :
\beginpicture
\setcoordinatesystem units <0.33truein, 0.33truein> point at 2.2 0.9
\setplotsymbol ({$\cal G$})
\plotsymbolspacing=9truept
\circulararc 315 degrees from 0 1 center at 0 0
\setplotsymbol ({$\cal T$})
\circulararc 315 degrees from 1 -1 center at 1 0
\endpicture
%   end of logo
%
\break
{\small\ifx\thesamplenumber\relax % sample?  
Volume \else Sample
\fi\thevolumenumber\ (\thevolumeyear)
\startpage--\finishpage\nl
Published: \publishdate}
\vglue 0.5truein plus 0.4fil minus 0.1truein

% title
{\parskip=0pt\leftskip 0pt plus 1fil\def\\{\par\smallskip}{\ifplaintex\large
\else\Large\fi\bf\thetitle}\par\medskip}   

\vglue 0pt plus 0.1fil 

% authors
{\parskip=0pt\leftskip 0pt plus 1fil\def\\{\par}{\sc\theauthors}
\par\medskip}

\vglue 0pt plus 0.1fil 

%address(es)
{\small\parskip=0pt\let\newline\\
{\leftskip 0pt plus 1fil\def\\{\par}{\sl\theaddress}\par}
\expandafter\ifx\theemail\relax    % email address?
\relax\else\vglue 5pt plus 0.02fil minus 2pt\def\\{\stdspace{\rm 
and}\stdspace} 
\cl{Email:\stdspace\tt\theemail}\fi
\ifx\theurl\relax                  % URL given?
\relax\else\vglue 5pt plus 0.02fil minus 2pt\def\\{\stdspace{\rm 
and}\stdspace}
\cl{URL:\stdspace\tt\theurl}\fi\par}

\vglue 7pt plus 0.3fil minus 3pt

{\bf Abstract}
\vglue 5pt plus 0.1fil minus 2pt

\theabstract

\vglue 7pt plus 0.3fil minus 3pt

{\bf AMS Classification numbers}\quad Primary:\quad \theprimaryclass

Secondary:\quad \thesecondaryclass

\vglue 5pt plus 0.3fil minus 2pt

{\bf Keywords}\quad \thekeywords

\vglue 10pt plus 0.5fil minus 5pt

{\small  Proposed: \theproposer\hfill Received: \receiveddate\nl
Seconded: \theseconders\hfill 
\ifx\reviseddate\relax                         % paper revised?
Accepted: \accepteddate                        % no
\else
Revised: \reviseddate                          % yes
\fi}
\eject
}       %  end of definition of \maketitlep

\let\maketitlepage\maketitlep
\let\maketitle\maketitlepage

%%% Headers and footers

\font\phead=cmsl9 scaled 950
\font\lhead=cmsl9 scaled 1050
\font\pnum=cmbx10 scaled 913
\font\lnum=cmbx10 
\font\pfoot=cmsl9 scaled 950
\font\lfoot=cmsl9 scaled 1050
\ifplaintex
\headline{\vbox to 0pt{\vskip -4.5mm\line{\small\phead\ifnum
\count0=\startpage ISSN 1364-0380 (on line)
1465-3060 (printed) \hfill {\pnum\folio}\else\ifodd\count0\def\\{ }% 
\ifx\theshorttitle\relax\thetitle\else\theshorttitle\fi\hfill{\pnum\folio}
\else\def\\{ and }{\pnum\folio}\hfill\ifx\theshortauthors\relax\theauthors
\else\theshortauthors\fi\fi\fi}\vss}}
\footline{\vbox to 0pt{\vglue 0mm\line{\small\pfoot\ifnum\count0=\startpage
\copyright\ \gtp\hfill\else
\gt, Volume \thevolumenumber\ (\thevolumeyear)\hfill\fi}\vss
}}
\else
\makeatletter
\def\@oddhead{{\small\lhead\ifnum\count0=\startpage ISSN 1364-0380 (on line)
1465-3060 (printed) \hfill {\lnum\number\count0}\else\ifodd\count0
\def\\{ }\ifx\theshorttitle\relax \thetitle \else\theshorttitle\fi\hfill
{\lnum\number\count0}\else\def\\{ and }{\lnum\number\count0}
\hfill\ifx\theshortauthors\relax 
\theauthors\else\theshortauthors\fi\fi\fi}}\def\@evenhead{\@oddhead}
\def\@oddfoot{\small\lfoot\ifnum\count0=\startpage\copyright\ \gtp\hfill\else
\gt, Volume \thevolumenumber\ (\thevolumeyear)\hfill\fi}
\def\@evenfoot{\@oddfoot}
\makeatother
\fi

   %%%comment out to create xxx header file

\newwrite\gtoutfile
\long\gdef\makeheadfile{  %%% start of definition of \makeheadfile
{\def\\{, }\def\s{ }
\immediate\openout\gtoutfile head.xxx
\immediate\write\gtoutfile{To: math@arxiv.org}
\immediate\write\gtoutfile{Subject: put or rep NNNNN:pppp}
\immediate\write\gtoutfile{--text follows this line--}
\immediate\write\gtoutfile{Proxy-for: \ifx\theasciiauthors\relax
\theauthors\else\theasciiauthors\fi\s<\ifx\theasciiemail\relax\theemail\else\theasciiemail\fi>}
\immediate\write\gtoutfile{\noexpand\\}
\immediate\write\gtoutfile{Authors: \ifx\theasciiauthors\relax
\theauthors\else\theasciiauthors\fi}
{\def\\{ }\immediate\write\gtoutfile{Title: \ifx\theasciititle\relax
\thetitle\else\theasciititle\fi}}
\immediate\write\gtoutfile{Subj-class: GT or SG or MG etc}
\immediate\write\gtoutfile{MSC-class: \theprimaryclass\ifx\thesecondaryclass\relax\else, \thesecondaryclass\fi}
\immediate\write\gtoutfile{Journal-ref: Geom. Topol. \thevolumenumber
(\thevolumeyear) \startpage-\finishpage}
\immediate\write\gtoutfile{Comments: Published by Geometry and Topology at}
\immediate\write\gtoutfile{\s\s http://www.maths.warwick.ac.uk/gt/GTVol\thevolumenumber/paper\thepapernumber.abs.html}
\immediate\write\gtoutfile{\noexpand\\}
\immediate\write\gtoutfile{}
\ifx\theasciiabstract\relax
\immediate\write\gtoutfile{\theabstract}\else
\immediate\write\gtoutfile{\theasciiabstract}\fi
\immediate\write\gtoutfile{}
\immediate\write\gtoutfile{\noexpand\\}
\immediate\write\gtoutfile{}
\immediate\closeout\gtoutfile}}  %%% end of definition of \makeheadfile

\def\maketitlepage{\maketitlep\makeheadfile}
\let\maketitle\maketitlepage

\lognumber{221}
\volumenumber{7}\papernumber{4}\volumeyear{2003}
\pagenumbers{155}{184}
\received{30 November 2001}
\revised{7 February 2003}
\accepted{13 March 2003}
\published{14 March 2003}
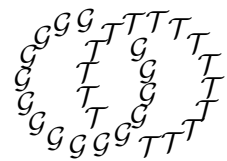
\proposed{Haynes Miller}
\seconded{Gunnar Carlsson, Thomas Goodwillie}

\input amsnames
\input amstex
\loadeusm
\input xy
\xyoption{all}
\let\cal\Cal      
\catcode`\@=12    
\let\\\par
\def\topmatter{\relax}
\def\endtopmatter{\maketitlepage}
\let\gttitle\title
\def\title#1\endtitle{\gttitle{#1}}
\let\gtauthor\author
\def\author#1\endauthor{\gtauthor{#1}}
\let\gtaddress\address
\def\address#1\endaddress{\gtaddress{#1}}
\let\gtemail\email
\def\email#1\endemail{\gtemail{#1}}
\def\subjclass#1\endsubjclass{\primaryclass{#1}}
\let\gtkeywords\keywords
\def\keywords#1\endkeywords{\gtkeywords{#1}}
\def\heading#1\endheading{{\def\S##1{\relax}\def\\{\relax\ignorespaces}
    \section{#1}}}
\def\head#1\endhead{\heading#1\endheading}

\def\subhead#1\endsubhead{\sh{#1}}
\def\subsubhead#1\endsubsubhead{\sh{#1}}
\def\specialhead#1\endspecialhead{\sh{#1}}
\def\demo#1{\rk{#1}\ignorespaces}
\def\enddemo{\ppar}
\let\remark\demo
\def\endremark{}

\let\example\demo
\def\endexample{\ppar}
\def\qed{\ifmmode\quad\sq\else\hbox{}\hfill$\sq$\par\goodbreak\rm\fi}  
\def\proclaim#1{\rk{#1}\sl\ignorespaces}
\def\endproclaim{\rm\ppar}
\def\cite#1{[#1]}
\newcount\itemnumber

\let\itemold\item
\def\item{\itemold{{\rm(\number\itemnumber)}}%
\global\advance\itemnumber by 1\ignorespaces}
\def\S{section~\ignorespaces}  %%  expand \S to "section"
\def\date#1\enddate{\relax}
\def\thanks#1\endthanks{\relax}   %%%  Move acknowledgements "manually"
\def\dedicatory#1\enddedicatory{\relax}  %%% to the end of the intro.
\def\rom#1{{\rm #1}}  % in some versions of amstex but not all
\let\footnote\plainfootnote

\def\Refs{\ppar{\large\bf References}\ppar\bgroup\leftskip=25pt
\frenchspacing\parskip=3pt plus2pt\small}       
\def\endRefs{\egroup}
\def\widestnumber#1#2{\relax}
\def\endrefitem{}
\def\refdef#1#2#3{\def#1{\leavevmode\unskip\endrefitem#2\def\endrefitem{#3}}}
\def\ref{\par}
\def\endref{\endrefitem\par\def\endrefitem{}}
\refdef\key{\noindent\llap\bgroup[}{]\ \ \egroup}
\refdef\no{\noindent\llap\bgroup[}{]\ \ \egroup}
\refdef\by{\bf}{\rm, }
\refdef\manyby{\bf}{\rm, }
\refdef\paper{\it}{\rm, }
\refdef\book{\it}{\rm, }
\refdef\jour{}{ }
\refdef\vol{}{ }
\refdef\yr{$(}{)$ }
\refdef\ed{(}{ Editor) }
\refdef\publ{}{ }
\refdef\inbook{from: ``}{'', }
\refdef\pages{}{ }
\refdef\page{}{ }
\refdef\paperinfo{}{ }
\refdef\bookinfo{}{ }
\refdef\publaddr{}{ }
\refdef\eds{(}{ Editors)}
\refdef\bysame{\hbox to 3 em{\hrulefill}\thinspace,}{ }
\refdef\toappear{(to appear)}{ }
\refdef\issue{no.\ }{ }
\newcount\refnumber\refnumber=1
\def\refkey#1{\expandafter\xdef\csname cite#1\endcsname{\number\refnumber}%
\global\advance\refnumber by 1}
\def\cite#1{[\csname cite#1\endcsname]}
\def\Cite#1{\csname cite#1\endcsname}  %% unbracketed \cite 
\def\key#1{\noindent\llap{[\csname cite#1\endcsname]\ \ }}

\refkey {Ad}
\refkey {AP}
\refkey {Be}
\refkey {BK}
\refkey {BHM}
\refkey {BM1}
\refkey {BM2}
\refkey {Bo}
\refkey {Du}
\refkey {Dw}
\refkey {DF1}
\refkey {DF2}
\refkey {DM}
\refkey {FH}
\refkey {HM}
\refkey {Ho}
\refkey {Ig}
\refkey {KR}
\refkey {Kn}
\refkey {Li1}
\refkey {Li2}
\refkey {MM}
\refkey {MS}
\refkey {MST}
\refkey {MT}
\refkey{MQR}
\refkey{Mil1}
\refkey{Mil2}
\refkey{Mis}
\refkey{Mit}
\refkey{Mo}
\refkey{Mu}
\refkey{Qu1}
\refkey{Qu2}
\refkey{Qu3}
\refkey{Qu4}
\refkey{Ra}
\refkey{Ro}
\refkey{RW}
\refkey{St}
\refkey{SV}
\refkey{Ta}
\refkey{To}
\refkey{Vo1}
\refkey{Vo2}
\refkey{Wa1}
\refkey{Wa2}
\refkey{Wa3}
\refkey{Wa4}
\refkey{Wa5}
\refkey{Wa6}

\define\Br{\operatorname{Br}}
\define\Cok{\operatorname{Cok}}
\define\C{\Bbb C}

\define\F{\Bbb F}
\define\Hom{\operatorname{Hom}}
\define\In{\operatorname{in}}
\define\Pic{\operatorname{Pic}}
\define\Q{\Bbb Q}
\define\Z{\Bbb Z}
\define\cok{\operatorname{cok}}
\define\et{\text{{\'e}t}}
\define\ev{\operatorname{ev}}
\define\hocolim{\operatorname{hocolim}}
\define\hofib{\operatorname{hofib}}
\define\pr{\operatorname{pr}}
\define\red{\operatorname{red}}
\define\tors{\operatorname{tors}}
\define\trc{\operatorname{trc}}
\define\trf{\operatorname{trf}}
\define\tr{\operatorname{tr}}
\redefine\H{\Bbb H}

\hyphenation{homo-logy}

\topmatter
\title The smooth Whitehead spectrum\\of a point at odd regular primes
\endtitle
\shorttitle{The smooth Whitehead spectrum}
\author John Rognes \endauthor
\address Department of Mathematics, University of Oslo\\N--0316 Oslo, Norway
\endaddress
\email rognes@math.uio.no \endemail
\abstract
Let $p$ be an odd regular prime, and assume that the Lichtenbaum--Quillen
conjecture holds for $K(\Z[1/p])$ at~$p$.  Then the $p$-primary
homotopy type of the smooth Whitehead spectrum $Wh(*)$ is described.
A suspended copy of the cokernel-of-J spectrum splits off, and the torsion
homotopy of the remainder equals the torsion homotopy of the fiber
of the restricted $S^1$-transfer map $t \: \Sigma \C P^\infty \to S$.
The homotopy groups of $Wh(*)$ are determined in a range of degrees, and the
cohomology of $Wh(*)$ is expressed as an $A$-module in all degrees, up to
an extension.  These results have geometric topological interpretations,
in terms of spaces of concordances or diffeomorphisms of highly connected,
high dimensional compact smooth manifolds.
\endabstract
\asciiabstract{%
Let p be an odd regular prime, and assume that the Lichtenbaum-Quillen
conjecture holds for K(Z[1/p])  at p.  Then the p-primary
homotopy type of the smooth Whitehead spectrum Wh(*) is described.
A suspended copy of the cokernel-of-J spectrum splits off, and the torsion
homotopy of the remainder equals the torsion homotopy of the fiber
of the restricted S^1-transfer map t: SigmaCP^infty--> S.
The homotopy groups of Wh(*) are determined in a range of degrees, and the
cohomology of Wh(*) is expressed as an A-module in all degrees, up to
an extension.  These results have geometric topological interpretations,
in terms of spaces of concordances or diffeomorphisms of highly connected,
high dimensional compact smooth manifolds.}
\primaryclass{19D10}
\secondaryclass{19F27, 55P42, 55Q52, 57R50, 57R80}
\keywords
Algebraic $K$-theory, topological cyclic homology, Lichtenbaum--Quillen
conjecture, transfer, $h$-cobordism, concordance, pseudoisotopy
\endkeywords
\asciikeywords{Algebraic K-theory, topological cyclic homology, 
Lichtenbaum-Quillen conjecture, transfer, h-cobordism, concordance,
pseudoisotopy}
\endtopmatter

\catcode`\@=\active

\document
\head Introduction \endhead

In this paper we study the smooth Whitehead spectrum $Wh(*)$ of
a point at an odd regular prime $p$, under the assumption that the
Lichtenbaum--Quillen conjecture for $K(\Z[1/p])$ holds at $p$.  This is
a reasonable assumption in view of recent work by Rost and Voevodsky.
The results admit geometric topological interpretations in terms of
the spaces of concordances (= pseudo-isotopies), $h$-cobordisms and
diffeomorphisms of high-dimensional compact smooth manifolds that are
as highly connected as their concordance stable range.  Examples of such
manifolds include discs and spheres.

Here is a summary of the paper.

We begin in \S2 by recalling Waldhausen's algebraic $K$-theory of
spaces \cite{Wa4}, Quillen's algebraic $K$-theory of rings \cite{Qu1},
the Lichtenbaum--Quillen conjecture in the strong formulation of Dwyer
and Friedlander \cite{DF1}, and a theorem of Dundas \cite{Du} about the
relative properties of the cyclotomic trace map to the topological cyclic
homology of B{\"o}kstedt, Hsiang and Madsen \cite{BHM}.

From \S3 and onwards we assume that $p$ is an odd regular prime and
that the Lichtenbaum--Quillen conjecture holds for $K(\Z[1/p])$ at $p$.
In~3.1 and~3.3 we then call on Tate--Poitou duality for {\'e}tale cohomology
\cite{Ta} to obtain a cofiber sequence
$$
j \vee \Sigma^{-2} ko @>>> Wh(*) @>\widetilde{\trc}>> \Sigma\C
P^\infty_{-1} @>>> \Sigma j \vee \Sigma^{-1} ko
\tag 1.1
$$
of implicitly $p$-completed spectra.  Here $\C P^\infty_{-1} =
Th(-\gamma^1)$ is a stunted complex projective
spectrum with one cell in each even dimension $\ge -2$, $j$ is the
connective image-of-J spectrum at $p$, and $ko$ is the connective real
$K$-theory spectrum.
In~3.6 we use this to obtain a splitting
$$
Wh(*) \simeq \Sigma c \vee (Wh(*)/\Sigma c)
\tag 1.2
$$
of the suspended cokernel-of-J spectrum $\Sigma c$ off from $Wh(*)$,
and in~3.8 we obtain a cofiber sequence
$$
\Sigma^2 ko @>>> Wh(*)/\Sigma c @>\tau>>
P_0\Sigma\overline{\C P}^\infty_{-1} @>>> \Sigma^3 ko \,,
\tag 1.3
$$
where $\pi_*(\tau)$ identifies the $p$-torsion in the homotopy of
$Wh(*)/\Sigma c$ with that of $\Sigma\overline{\C P}^\infty_{-1}$.
The latter spectrum equals the homotopy fiber of the restricted
$S^1$-transfer map
$$
t \: \Sigma \C P^\infty \to S \,.
$$
Hence the homotopy of $Wh(*)$ is as complicated as the (stable) homotopy
of infinite complex projective space $\C P^\infty$, and the associated
transfer map above.

In \S4 we make a basic homotopical analysis, following Mosher \cite{Mo}
and Knapp \cite{Kn}, to compute $\pi_* \Sigma\overline{\C P}^\infty_{-1}$
and thus $\pi_* Wh(*)$ at $p$ in degrees up to $|\beta_2|-2 = (2p+1)q-4$,
where $q=2p-2$ as usual.  See~4.7 and~4.9.  The first $p$-torsion to
appear in $\pi_m Wh(*)$ is $\Z/p$ for $m = 4p-2$ when $p\ge5$, and
$\Z/3\{\Sigma\beta_1\}$ for $m=11$ when $p=3$.

In \S5 we make the corresponding mod~$p$ cohomological analysis and
determine $H^*(Wh(*); \F_p)$ as a module over the Steenrod algebra is
all degrees, up to an extension.  See~5.4 and~5.5.  The extension is
trivial for $p=3$, and nontrivial for $p\ge5$.  Taken together, this
homotopical and cohomological information gives a detailed picture of
the homotopy type $Wh(*)$.

In \S6 we recall the relation between the Whitehead spectrum $Wh(*)$,
the concordance space $C(M)$ and the diffeomorphism group
$DIFF(M)$ of suitably highly connected and high dimensional compact
smooth manifolds~$M$.  As a sample application we show in~6.3 that
for $p\ge5$ and $M$ a compact smooth $k$-connected $n$-manifold with
$k \ge 4p-2$ and $n \ge 12p-5$, the first $p$-torsion in the homotopy
of the smooth concordance space $C(M)$ is $\pi_{4p-4} C(M)_{(p)} \cong
\Z/p$.  Specializing to $M = D^n$ we conclude in~6.4 that $\pi_{4p-4}
DIFF(D^{n+1})$ or $\pi_{4p-4} DIFF(D^n)$ contains an element of order
exactly~$p$.  Comparable results hold for $p=3$.

A 2-primary analog of this study was presented in \cite{Ro}.  Related
results on the homotopy fiber of the linearization map $L \: A(*) \to
K(\Z)$ were given in \cite{KR}.

\head Algebraic $K$-theory and topological cyclic homology \endhead

\subhead Algebraic $K$-theory of spaces \endsubhead
Let $A(X)$ be Waldhausen's algebraic $K$-theory spectrum [\Cite{Wa4}, \S2.1]
of a space $X$.  There is a natural cofiber sequence [\Cite{Wa4}, \S3.3],
\cite{Wa5}
$$
\Sigma^\infty(X_+) @>\eta_X>> A(X) @>\pi>> Wh(X) \,,
$$
where $Wh(X) = Wh^{DIFF}(X)$ is the smooth Whitehead spectrum of $X$,
and a natural trace map \cite{Wa2} $\tr_X \: A(X) \to \Sigma^\infty(X_+)$
which splits the above cofiber sequence up to homotopy.  Let $\iota
\: Wh(X) \to A(X)$ be the corresponding homotopy section to $\pi$.
When $X=*$ is a point, $\Sigma^\infty(*_+) = S$ is the sphere spectrum,
and the splitting simplifies to $A(*) \simeq S \vee Wh(*)$.

\subhead Topological cyclic homology of spaces \endsubhead
Let $p$ be a prime and let $TC(X; p)$ be B{\"o}kstedt, Hsiang and
Madsen's topological cyclic homology [\Cite{BHM}, 5.12(i)] of the space $X$.
There is a natural cofiber sequence [\Cite{BHM}, 5.17]
$$
\hofib(\trf_{S^1}) @>\iota>> TC(X; p) @>\beta_X>> \Sigma^\infty(\Lambda X_+)
$$
after $p$-adic completion, where $\Lambda X$ is the free loop space of
$X$ and
$$
\trf_{S^1} \: \Sigma^\infty(\Sigma (ES^1 \times_{S^1} \Lambda X)_+)
\to \Sigma^\infty(\Lambda X_+)
$$
is the dimension-shifting $S^1$-transfer map for the canonical
$S^1$-bundle $ES^1 \times \Lambda X \to ES^1 \times_{S^1} \Lambda X$;
see e.g.~[\Cite{MS}, \S2].  When $X=*$ the $S^1$-transfer map simplifies to
$\trf_{S^1} \: \Sigma^\infty \Sigma \C P^\infty_+ \to S$.  Its homotopy
fiber is $\Sigma \C P^\infty_{-1}$ [\Cite{MS}, \S3], where the stunted
complex projective spectrum $\C P^\infty_{-1} = Th(-\gamma^1 \downarrow
\C P^\infty)$ is defined as the Thom spectrum of minus the tautological
line bundle over $\C P^\infty$.  The map $\iota$ identifies $\Sigma \C
P^\infty_{-1}$ with the homotopy fiber of $\beta_* \: TC(*; p) \to S$,
after $p$-adic completion.

We can think of $\C P^\infty_{-1}$ as a CW spectrum, with $2k$-skeleton
$\C P^k_{-1} = Th(-\gamma^1 \downarrow \C P^{k+1})$.  By James periodicity
$\Sigma^{2n} \C P^k_{-1} \simeq \C P^{n+k}_{n-1} = \C P^{n+k}/\C P^{n-2}$
whenever $n$ is a multiple of a suitable natural number that depends
on $k$.  From this it follows that integrally $H_*(\C P^\infty_{-1})
\cong \Z\{b_k \mid k\ge-1\}$ and $H^*(\C P^\infty_{-1}) \cong \Z\{y^k
\mid k\ge-1\}$ with $y^k$ dual to $b_k$, both in degree~$2k$.  In mod~$p$
cohomology the Steenrod operations act by $P^i(y^k) = {k \choose i}
y^{k+(p-1)i}$ and $\beta(y^k) = 0$.  In particular $P^i(y^{-1}) =
(-1)^i y^{-1+(p-1)i} \ne 0$ for all $i\ge0$.

\subhead The cyclotomic trace map for spaces \endsubhead Let $\trc_X \:
A(X) \to TC(X; p)$ be the natural cyclotomic trace map of B{\"o}kstedt,
Hsiang and Madsen [\Cite{BHM}, 5.12(ii)].  It lifts the Waldhausen trace
map, in the sense that $\tr_X \simeq \ev \circ \beta_X \circ \trc_X$,
where $\ev \: \Sigma^\infty(\Lambda X_+) \to \Sigma^\infty(X_+)$
evaluates a free loop at a base point.  Hence there is a map of (split)
cofiber sequences of spectra:
$$
\xymatrix{
Wh(X) \ar[r]^\iota \ar[d]^{\widetilde{\trc}} &
A(X) \ar[r]^-{\tr_X} \ar[d]^{\trc_X} &
\Sigma^\infty(X_+) \ar[d]^{=} \\
{}\hofib(\ev \circ \beta_X) \ar[r]^\iota &
TC(X; p) \ar[r]^-{\ev \circ \beta_X} &
\Sigma^\infty(X_+)
}
$$
after $p$-adic completion.  When $X=*$ the left hand square simplifies
as follows:

\proclaim{Theorem 2.1}\qua{\rm (Waldhausen, B{\"o}kstedt--Hsiang--Madsen)}\qua
There is a homotopy Cartesian square
$$
\xymatrix{
Wh(*) \ar[r]^\iota \ar[d]^{\widetilde{\trc}} &
A(*) \ar[d]^{\trc_*} \\
\Sigma \C P^\infty_{-1} \ar[r]^\iota &
TC(*; p)
}
$$
after $p$-adic completion.  Hence there is a $p$-complete equivalence
$\hofib(\widetilde{trc}) \simeq \hofib(\trc_*)$.
\endproclaim

\subhead Algebraic $K$-theory of rings \endsubhead
Let $K(R)$ be Quillen's algebraic $K$-theory spectrum of a ring $R$
[\Cite{Qu1}, \S2].  When $R$ is commutative, Noetherian and $1/p \in R$
the {\'e}tale $K$-theory spectrum $K^{\et}(R)$ of Dwyer and Friedlander
[\Cite{DF1}, \S4] is defined, and comes equipped with a natural comparison
map $\phi \: K(R) \to K^{\et}(R)$.  By construction $K^{\et}(R)$ is a
$p$-adically complete $K$-local spectrum \cite{Bo}.  Let $R$ be the ring
of $p$-integers in a local or a global field of characteristic $\ne p$.
The Lichtenbaum--Quillen conjecture \cite{Li1}, \cite{Li2}, \cite{Qu3}
for $K(R)$ at $p$, in the strong form due to Dwyer and Friedlander,
then asserts:

\proclaim{Conjecture 2.2}\qua {\rm (Lichtenbaum--Quillen)}\qua
The comparison map $\phi$ induces a homotopy equivalence
$$
P_1 \phi^\wedge_p \: P_1 K(R)^\wedge_p @>>> P_1 K^{\et}(R)
$$
of $0$-connected covers after $p$-adic completion.
\endproclaim

Here $P_n E$ denotes the $(n-1)$-connected cover of any spectrum $E$.
In the cases of concern to us the $p$-completed map $\phi^\wedge_p$
will also induce an isomorphism in degree~0, so the covers $P_1$ above
can be replaced by $P_0$.

The conjecture above has been proven for $p=2$ by Rognes and Weibel
[\Cite{RW}, 0.6], based on Voevodsky's proof \cite{Vo1}, \cite{Vo2} of the
Milnor conjecture.  The odd-primary version of this conjecture would
follow \cite{SV} from results on the Bloch--Kato conjecture \cite{BK}
announced as ``in preparation'' by Rost and Voevodsky, but have not yet
formally appeared.

\subhead Topological cyclic homology of rings \endsubhead
Let $TC(R; p)$ be B{\"o}kstedt, Hsiang and Madsen's topological cyclic
homology of a (general) ring $R$.  There is a natural cyclotomic trace map
$\trc_R \: K(R) \to TC(R; p)$.  When $X$ is a based connected space with
fundamental group $\pi = \pi_1(X)$, and $R = \Z[\pi]$ is the group ring,
there are natural linearization maps $L \: A(X) \to K(R)$ [\Cite{Wa1},
\S2] and $L \: TC(X; p) \to TC(R; p)$ which commute with the cyclotomic
trace maps.  Moreover, by Dundas \cite{Du} the square
$$
\xymatrix{
A(X) \ar[r]^L \ar[d]^{\trc_X} &
K(R) \ar[d]^{\trc_R} \\
TC(X; p) \ar[r]^L &
TC(R; p)
}
$$
is homotopy Cartesian after $p$-adic completion.
In the special case when $X = *$ and $R = \Z$ this simplifies to:

\proclaim{Theorem 2.3}\qua {\rm (Dundas)}\qua
There is a homotopy Cartesian square
$$
\xymatrix{
A(*) \ar[r]^L \ar[d]^{\trc_*} &
K(\Z) \ar[d]^{\trc_{\Z}} \\
TC(*; p) \ar[r]^L &
TC(\Z; p)
}
$$
after $p$-adic completion.  Hence there is a $p$-complete equivalence
$\hofib(\trc_*) \simeq \hofib(\trc_{\Z})$.
\endproclaim

\subhead The cyclotomic trace map for rings \endsubhead
When $k$ is a perfect field of characteristic $p>0$, $W(k)$ its ring
of Witt vectors, and $R$ is an algebra of finite rank over $W(k)$, then
by Hesselholt and Madsen [\Cite{HM}, Thm.~D] there is a cofiber sequence
of spectra
$$
K(R) @>\trc_R>> TC(R; p) @>>> \Sigma^{-1} HW(R)_F
$$
after $p$-adic completion.  Here $W(R)_F$ equals the coinvariants of
the Frobenius action on the Witt ring of $R$, and $\Sigma^{-1} HW(R)_F$
is the associated desuspended Eilenberg--Mac\,Lane spectrum.  The Witt
ring of $k = \F_p$ is the ring $W(\F_p) = \Z_p$ of $p$-adic integers,
so the above applies to $R = \Z_p[\pi]$ for finite groups $\pi$.
In particular, when $X = *$ and $\pi = 1$ there is a cofiber sequence
$$
K(\Z_p) @>\trc_{\Z_p}>> TC(\Z_p; p) @>>> \Sigma^{-1} H\Z_p
$$
after $p$-adic completion.  This uses that $W(\Z_p)_F \cong \Z_p$.

\subhead The completion map \endsubhead
Let $\kappa \: \Z \to \Z_p$ and $\kappa' \: \Z[1/p] \to \Q_p$ be the
$p$-completion homomorphisms, where $\Q_p$ is the field of $p$-adic
numbers.  By naturality of $\trc_R$ with respect to $\kappa$ there is
a commutative square
$$
\xymatrix{
K(\Z) \ar[r]^\kappa \ar[d]^{\trc_{\Z}} &
K(\Z_p) \ar[d]^{\trc_{\Z_p}} \\
TC(\Z; p) \ar[r]^-{\simeq} &
TC(\Z_p; p) \,.
}
$$
The lower map is a $p$-adic equivalence, since topological cyclic
homology is insensitive to $p$-adic completion, cf. [\Cite{HM}, \S6].
Hence there is a cofiber sequence of homotopy fibers
$$
\hofib(\kappa) @>>> \hofib(\trc_{\Z}) @>>> \Sigma^{-2} H\Z_p \,.
$$
By the localization sequences in $K$-theory [\Cite{Qu1}, \S5] there is a
homotopy Cartesian square
$$
\xymatrix{
K(\Z) \ar[r] \ar[d]^\kappa & K(\Z[1/p]) \ar[d]^{\kappa'} \\
K(\Z_p) \ar[r] & K(\Q_p)
}
$$
so $\hofib(\kappa) \simeq \hofib(\kappa')$.

\subhead Topological $K$-theory and related spectra \endsubhead
Let $ko$ and $ku$ be the connective real and complex topological $K$-theory
spectra, respectively.  There is a complexification map $c \: ko \to ku$,
and a cofiber sequence
$$
\Sigma ko @>\eta>> ko @>c>> ku @>r\beta^{-1}>> \Sigma^2 ko
$$
related to real Bott periodicity, cf.~[\Cite{MQR}, V.5.15].  Here $\eta$
is multiplication by the stable Hopf map $\eta \: S^1 \to S^0$, which
is null-homotopic at odd primes, $\beta \: \Sigma^2 ku \to ku$ covers
the Bott equivalence, and $r \: ku \to ko$ is realification.

Suppose $p$ is odd, and let $q = 2p-2$. There are splittings
$ku_{(p)} \simeq \bigvee_{i=0}^{p-2} \Sigma^{2i} \ell$ and
$$
ko_{(p)} \simeq \bigvee_{i=0}^{(p-3)/2} \Sigma^{4i} \ell \,,
\tag 2.4
$$
where $\ell$ is the connective $p$-local Adams summand of $ku$ \cite{Ad}.
There is a cofiber sequence $\Sigma^q \ell \to \ell \to H\Z_{(p)}$ that
identifies $\Sigma^q \ell$ with $P_q \ell$.  Let $r$ be a topological
generator of the $p$-adic units $\Z_p^*$, and let $\psi^r$ be the
Adams operation.  The $p$-local image-of-J spectrum $j$ is defined [\Cite{MQR},
V.5.16] by the cofiber sequence
$$
j @>>> \ell @>\psi^r-1>> \Sigma^q \ell \,.
$$
We now briefly write $S$ for the $p$-local sphere spectrum.  There is
a unit map $e \: S \to j$ representing (minus) the Adams $e$-invariant
on homotopy \cite{Qu4}, and the $p$-local cokernel-of-J spectrum $c$
is defined by the cofiber sequence
$$
c @>f>> S @>e>> j \,.
\tag 2.5
$$
Here $e$ induces a split surjection on homotopy, so $\pi_*(f)$ is split
injective.  The map $e$ identifies $j$ with the connective cover $P_0
L_K S$ of the $K$-localization of $S$, localized at $p$ [\Cite{Bo}, 4.3].

\proclaim{Lemma 2.6}
Suppose that $n \le 2q$.  If $n \ne q+1$ there are no essential spectrum
maps $H\Z_{(p)} \to \Sigma^n \ell$.  If $n = q+1$ the group of spectrum
maps $H\Z_{(p)} \to \Sigma^{q+1} \ell$ is $\Z_{(p)}$, generated by the
connecting map $\partial$ of the cofiber sequence $\Sigma^q \ell \to
\ell \to H\Z_{(p)}$.
\endproclaim

\proclaim{Lemma 2.7}
There are no essential spectrum maps $\Sigma^n \ell \to j$ for $n\ge0$
even.  Hence there are no essential spectrum maps $\Sigma ko_{(p)}
\to \Sigma j$.
\endproclaim

The proofs are easy, using \cite{Mis} for~2.6, and [\Cite{MST}, Cor.~C]
or [\Cite{Mit}, 2.4] for~2.7.

\head Splittings at odd regular primes \endhead

\subhead The completion map in {\'e}tale $K$-theory \endsubhead
When $R = \Z[1/p]$ and $p$ is an odd regular prime there is a homotopy
equivalence $P_0 K^{\et}(\Z[1/p]) \simeq j \vee \Sigma ko$ after
$p$-adic completion [\Cite{DF2}, 2.3].  Taking into account that $\phi$
is an equivalence in degree~0 and that $K(\Z[1/p])$ has finite type
\cite{Qu2}, the Lichtenbaum--Quillen conjecture for $\Z[1/p]$ at $p$
amounts to the assertion that $K(\Z[1/p]) \simeq j \vee \Sigma ko$
after $p$-localization.  By the localization sequence in $K$-theory, this
is equivalent to the assertion that $K(\Z) \simeq j \vee \Sigma^5 ko$,
after $p$-localization.

Hereafter we (often implicitly) {\bf complete all spectra} at~$p$.

When $R = \Q_p$ and $p$ is an odd prime there is a $p$-adic
equivalence $P_0 K^{\et}(\Q_p) \simeq j \vee \Sigma j \vee \Sigma ku$.
The Lichtenbaum--Quillen conjecture for $\Q_p$ at $p$ asserts that
$K(\Q_p) \simeq j \vee \Sigma j \vee \Sigma ku$ [\Cite{DM}, 13.3], which
again is equivalent to the assertion that $K(\Z_p) \simeq j \vee \Sigma
j \vee \Sigma^3 ku$, after $p$-adic completion.  This is now a theorem,
following from the calculation by B{\"o}kstedt and Madsen of $TC(\Z;
p)$ [\Cite{BM1}, 9.17], \cite{BM2}.

\proclaim{Proposition 3.1}
Let $p$ be an odd regular prime.  There are $p$-adic equivalences $P_0
K^{\et}(\Z[1/p]) \simeq j \vee \Sigma ko$ and $P_0 K^{\et}(\Q_p) \simeq
j \vee \Sigma j \vee \Sigma ku$ such that
$$
\kappa' \: P_0 K^{\et}(\Z[1/p]) @>>> P_0 K^{\et}(\Q_p)
$$
is homotopic to the wedge sum of the identity $id \: j \to j$, the zero
map $* \to \Sigma j$, and the suspended complexification map $\Sigma c \:
\Sigma ko \to \Sigma ku$.  Thus $\hofib(\kappa') \simeq j \vee \Sigma^2 ko$.
\endproclaim

\demo{Proof}
Taking the topological generator $r$ to be a prime power, there is a
reduction map $\red \: P_0 K^{\et}(\Q_p) \to K(\F_r) \simeq j$ after
$p$-adic completion [\Cite{DM}, \S13], such that the composite map
$$
S @>\eta>> K(\Z[1/p]) @>\phi>> P_0 K^{\et}(\Z[1/p]) @>\kappa'>> P_0
K^{\et}(\Q_p) @>\red>> j
$$
is homotopic to $e$.  Since $K^{\et}(\Z[1/p])$ is $K$-local, $\phi\eta$
also factors through $e$.  These maps split off a common copy of
$j$ from $P_0 K^{\et}(\Z[1/p])$ and $P_0 K^{\et}(\Q_p)$.  There are no
essential spectrum maps $\Sigma ko \to \Sigma j$ by~2.7, so after $p$-adic
completion $\kappa'$ is homotopic to a wedge sum of maps $id \: j \to
j$, $* \to \Sigma j$ and a map $\kappa'' \: \Sigma ko \to \Sigma ku$.
Any such $\kappa''$ lifts over $\Sigma c \: \Sigma ko \to \Sigma ku$, so
it suffices to show that $\pi_{2i-1}(\kappa'')$ is a $p$-adic isomorphism
for all odd $i\ge1$.

Equivalently we must show that $\kappa'$ induces an isomorphism on homotopy
modulo torsion subgroups in degree $2i-1$ for all odd $i>1$, or that
$$
K^{\et}_{2i-1}(\kappa'; \Q_p/\Z_p) \: K^{\et}_{2i-1}(\Z[1/p]; \Q_p/\Z_p)
@>>> K^{\et}_{2i-1}(\Q_p; \Q_p/\Z_p)
$$
is injective.  This equals the completion map
$$
\kappa' \: H_{\et}^1(\Z[1/p]; \Q_p/\Z_p(i)) \to H_{\et}^1(\Q_p; \Q_p/\Z_p(i))
$$
in {\'e}tale cohomology, by the collapsing spectral sequence in
[\Cite{DF}, 5.1].  By the 9-term exact sequence expressing Tate--Poitou
duality [\Cite{Ta}, 3.1], [\Cite{Mil2}, I.4.10], its kernel is a quotient
of $A^\# = H_{\et}^2(\Z[1/p]; \Z_p(1-i))^\#$, where $A^\# = \Hom(A,
\Q/\Z)$ denotes the Pontryagin dual of an abelian group $A$.  But $A =
H_{\et}^2(\Z[1/p]; \Z_p(1-i))$ is an abelian pro-$p$-group, with $A/p
\cong H_{\et}^2(\Z[1/p]; \Z/p(1-i))$ contained as a direct summand in $B =
H_{\et}^2(\Z[1/p, \zeta_p]; \Z/p)$, which is independent of $i$.  Here $R
= \Z[1/p, \zeta_p]$ is the ring of $p$-integers in the $p$-th cyclotomic
field $\Q(\zeta_p)$.  Kummer theory gives a short exact sequence
$$
0 \to \Pic(R)/p @>>> B @>>> \{p\} \Br(R) \to 0
$$
where $\Pic(R)$ and $\Br(R)$ are the Picard and Brauer groups of
$R$, respectively.  (See [\Cite{Mil2}, \S IV] and \cite{Ho}.)  Here $\Pic(R)/p
= 0$ because $p$ is a regular prime, and $\{p\} \Br(R) = \ker(p \: \Br(R)
\to \Br(R)) = 0$ because $p$ is odd and $(p)$ does not split in $R$
[\Cite{Mil1}, p.~109], so $B = 0$.  Thus $A/p = 0$ and it follows that $A =
0$, since $A$ is an abelian pro-$p$-group.
\qed
\enddemo

\subhead The fiber of the cyclotomic trace map \endsubhead
Hereafter we make the following standing assumption.

\proclaim{Hypothesis 3.2}
{\rm(a)}\qua $p$ is an odd regular prime, and

{\rm(b)}\qua the Lichtenbaum--Quillen conjecture~2.2 holds for $K(\Z[1/p])$ at~$p$.
\endproclaim

\proclaim{Proposition 3.3}
There is a homotopy equivalence $\hofib(\trc_{\Z}) \simeq j \vee
\Sigma^{-2} ko$ after $p$-adic completion.
\endproclaim

\demo{Proof}
By assumption $\kappa' \: K(\Z[1/p]) \to K(\Q_p)$ agrees with
$$
\kappa' \: P_0 K^{\et}(\Z[1/p]) \to P_0 K^{\et}(\Q_p)
$$
after $p$-adic completion, so we have a cofiber sequence
$$
j \vee \Sigma^2 ko @>>> \hofib(\trc_{\Z}) @>>> \Sigma^{-2} H\Z_p \,.
$$
The connecting map $\Sigma^{-2} H\Z_p \to \Sigma j \vee \Sigma^3 ko$ is
homotopic to a wedge sum of maps $\Sigma^{-2} H\Z_p \to \Sigma j$ and
$\Sigma^{-2} H\Z_p \to \Sigma^{4i-1} \ell$ for $1 \le i \le (p-1)/2$.
All such maps are null-homotopic by~2.6, with the exception of the map
$\partial' \: \Sigma^{-2} H\Z_p \to \Sigma^{2p-3} \ell$ corresponding to
$i = (p-1)/2$.

We claim that multiplication by $v_1$ acts nontrivially from degree~$-2$
to degree~$2p-4$ in $\pi_*(\hofib(\trc_{\Z}); \Z/p)$, from which it
follows that $\partial'$ is a $p$-adic unit times the connecting map
$\partial$ in the cofiber sequence $\Sigma^{q-2} \ell \to \Sigma^{-2}
\ell \to \Sigma^{-2} H\Z_p$.  This implies that
$$
\hofib(\trc_{\Z}) \simeq j \vee \Sigma^{-2} \ell \vee
\bigvee_{i=1}^{(p-3)/2} \Sigma^{4i-2} \ell \simeq j \vee \Sigma^{-2}
ko \,.
$$
To prove the claim, consider the homotopy Cartesian squares in~2.1 and~2.3.
In the Atiyah--Hirzebruch spectral sequence
$$
E^2_{s, t} = H_s(\C P^\infty_{-1}; \pi_t(S; \Z/p))
\Longrightarrow \pi_{s+t}(\C P^\infty_{-1}; \Z/p)
$$
there is a first differential $d^{2p-2}(b_{p-2}) = \alpha_1 b_{-1}$,
so we find $\pi_{-2}(\C P^\infty_{-1}; \Z/p) \cong \Z/p\{b_{-1}\}$
and $\pi_{2p-4}(\C P^\infty_{-1}; \Z/p) \cong \Z/p\{v_1b_{-1}\}$.
Hence multiplication by $v_1$ acts nontrivially from
$$
\pi_{-1}(TC(*; p); \Z/p) \cong \Z/p\{\Sigma b_{-1}\}
$$
to
$$
\pi_{2p-3}(TC(*; p); \Z/p) \cong \Z/p\{\alpha_1, \Sigma v_1b_{-1}\} \,,
$$
also modulo the image from the unit map $\eta \: S \to TC(*; p)$.

The map $L \: S \to H\Z$ is $(2p-3)$-connected, hence so is $L \: TC(*; p)
\to TC(\Z; p)$ by [\Cite{BM1}, 10.9] and \cite{Du}.  Here $\pi_{2p-3}(TC(\Z;
p); \Z/p) \cong \Z/p\{\alpha_1\} \oplus \Z/p$ since $P_0 TC(\Z; p) \simeq
K(\Z_p) \simeq j \vee \Sigma j \vee \Sigma^3 ku$.  So the surjection
$\pi_{2p-3}(L; \Z/p)$ is in fact a bijection, and multiplication by $v_1$
acts nontrivially from $\pi_{-1}(TC(\Z; p); \Z/p)$ to $\pi_{2p-3}(TC(\Z;
p); \Z/p)$, also modulo the image from the unit map $\eta \: S \to
TC(\Z; p)$.

By the assumed $p$-adic equivalence $K(\Z) \simeq j \vee \Sigma^5
ko$, this image equals the image from the cyclotomic trace map
$\trc_{\Z} \: K(\Z) \to TC(\Z; p)$.  Hence we can pass to cofibers,
and conclude that multiplication by $v_1$ acts nontrivially from
$\pi_{-2}(\hofib(\trc_{\Z}); \Z/p)$ to $\pi_{2p-4}(\hofib(\trc_{\Z});
\Z/p)$, as claimed.
\qed
\enddemo

We let $d$ be the homotopy cofiber map of $\widetilde{\trc}$.
Combining~2.1, 2.3 and~3.3 we have:

\proclaim{Corollary 3.4}
There is a diagram of horizontal cofiber sequences:
$$
\xymatrix{
j \vee \Sigma^{-2} ko \ar[r] \ar[d]^{=} &
Wh(*) \ar[r]^-{\widetilde{\trc}} \ar[d]^\iota &
\Sigma \C P^\infty_{-1} \ar[r]^-d \ar[d]^\iota &
\Sigma j \vee \Sigma^{-1} ko \ar[d]^{=} \\
j \vee \Sigma^{-2} ko \ar[r] \ar[d]^{=} &
A(*) \ar[r]^-{\trc_*} \ar[d]^L &
TC(*; p) \ar[r] \ar[d]^L &
\Sigma j \vee \Sigma^{-1} ko \ar[d]^{=} \\
j \vee \Sigma^{-2} ko \ar[r] &
K(\Z) \ar[r]^-{\trc_{\Z}} &
TC(\Z; p) \ar[r] &
\Sigma j \vee \Sigma^{-1} ko \,.
}
$$
\endproclaim

\subhead The restricted $S^1$-transfer map \endsubhead
There is a stable splitting $\In_1 \vee \In_2 \: S^1 \vee \Sigma \C
P^\infty \simeq \Sigma \C P^\infty_+$.  Let the restricted $S^1$-transfer
map $t = \trf_{S^1} \circ \In_2 \: \Sigma \C P^\infty \to S$ be the
restriction of $\trf_{S^1}$ to the second summand [\Cite{Mu}, \S2].
The restriction to the first summand is the stable Hopf map $\eta =
\trf_{S^1} \circ \In_1 \: S^1 \to S^0$, which is null-homotopic at odd
primes.  Hence the inclusion $\In_1$ lifts to a map $\Sigma b_0 \: S^1
\to \hofib(\trf_{S^1}) = \Sigma \C P^\infty_{-1}$, with Hurewicz image
$\Sigma b_0 \in H_1(\Sigma \C P^\infty_{-1})$.

Dually the projection $\pr_1 \: \Sigma \C P^\infty_+ \to S^1$ yields a
map $\Sigma y^0 \: \Sigma \C P^\infty_{-1} \to S^1$ with dual Hurewicz
image $\Sigma y^0 \in H^1(\Sigma \C P^\infty_{-1})$.  We obtain a diagram
of horizontal and split vertical cofiber sequences:
$$
\xymatrix{
S^1 \ar@<0.5ex>[d]^{\Sigma b_0} \ar[r]^{=} &
S^1 \ar@<0.5ex>[d]^{\In_1} \\
\Sigma \C P^\infty_{-1} \ar@<0.5ex>[d] \ar[r] \ar@<0.5ex>[u]^{\Sigma y^0} &
\Sigma \C P^\infty_+ \ar@<0.5ex>[d]^{\pr_2}
	\ar[r]^-{\trf_{S^1}} \ar@<0.5ex>[u]^{\pr_1} &
S \ar[d]^{=} \\
{}\hofib(t) \ar[r] \ar@<0.5ex>[u] &
\Sigma \C P^\infty \ar[r]^-t \ar@<0.5ex>[u]^{\In_2} &
S \,.
}
\tag 3.5
$$
Writing $\overline{\C P}^\infty_{-1}$ for the homotopy cofiber of $b_0 \:
S \to \C P^\infty_{-1}$, we have $\hofib(t) \simeq \Sigma \overline{\C
P}^\infty_{-1}$.  Then $H_*(\Sigma\overline{\C P}^\infty_{-1})
= \Z\{\Sigma b_k \mid k\ge-1, k\ne0\}$ and $H^*(\Sigma\overline{\C
P}^\infty_{-1}) \cong \Z\{\Sigma y^k \mid k\ge-1, k\ne0\}$.

It has been shown by Knapp \cite{Kn} that $\pi_*(t) \: \pi_*(\Sigma
\C P^\infty) \to \pi_*(S)$ is surjective for $0 < * < |\beta_{p+1}|
= p(p+2)q-2$, so the homotopy of $\Sigma\overline{\C P}^\infty_{-1}$
is as well understood in this range as that of $\Sigma\C P^\infty$.

\subhead The suspended cokernel-of-J spectrum \endsubhead
We can split off the suspension of the cofiber sequence~(2.5) defining
the cokernel-of-J from the top cofiber sequence in~3.4.

\proclaim{Proposition 3.6}
There is a diagram of horizontal and split vertical cofiber sequences:
$$
\xymatrix{
j \ar[r] \ar@<0.5ex>[d]^{\In_1} &
\Sigma c \ar[r]^{\Sigma f} \ar@<0.5ex>[d]^g &
S^1 \ar[r]^{\Sigma e} \ar@<0.5ex>[d]^{\Sigma b_0} &
\Sigma j \ar@<0.5ex>[d]^{\In_1} \\
j \vee \Sigma^{-2} ko \ar[r] \ar@<0.5ex>[d]^{\pr_2} \ar@<0.5ex>[u]^{\pr_1} &
Wh(*) \ar[r]^-{\widetilde{\trc}} \ar@<0.5ex>[d] \ar@<0.5ex>[u] &
\Sigma \C P^\infty_{-1} \ar[r]^-d \ar@<0.5ex>[d] \ar@<0.5ex>[u]^{\Sigma y^0} &
\Sigma j \vee \Sigma^{-1} ko \ar@<0.5ex>[d]^{\pr_2} \ar@<0.5ex>[u]^{\pr_1} \\
\Sigma^{-2} ko \ar[r] \ar@<0.5ex>[u]^{\In_2} &
Wh(*)/\Sigma c \ar[r] \ar@<0.5ex>[u] &
\Sigma\overline{\C P}^\infty_{-1} \ar[r] \ar@<0.5ex>[u] &
\Sigma^{-1} ko \ar@<0.5ex>[u]^{\In_2} \,.
}
$$
In particular there is a splitting
$$
Wh(*) \simeq \Sigma c \vee (Wh(*)/\Sigma c)
$$
where $Wh(*)/\Sigma c$ is defined as the homotopy cofiber of $g$.
\endproclaim

\demo{Proof}
The composite $d \circ \Sigma b_0$ represents the generator of
$\pi_1(\Sigma j \vee \Sigma^{-1} ko)$, hence factors as $\In_1 \circ
\Sigma e \: S^1 \to \Sigma j \to \Sigma j \vee \Sigma^{-1} ko$.  We define
$g \: \Sigma c \to Wh(*)$ as the induced map of homotopy fibers.  It is
well-defined up to homotopy since $\pi_2(\Sigma j \vee \Sigma^{-1}
ko) = 0$.  This explains the downward cofiber sequences of the diagram.

To split $g$ we must show that $\pr_1 \circ \, d$ factors as
$\Sigma e \circ \Sigma y^0$, or equivalently that the composite
$$
\Sigma\overline{\C P}^\infty_{-1} \to \Sigma \C P^\infty_{-1} @>d>>
\Sigma j \vee \Sigma^{-1} ko @>\pr_1>> \Sigma j
$$
is null-homotopic.  But this map lies in a zero group, because in the
Atiyah--Hirzebruch spectral sequence
$$
E^2_{s, t} = H^{-s}(\Sigma\overline{\C P}^\infty_{-1}; \pi_t(\Sigma j))
\Longrightarrow [\Sigma\overline{\C P}^\infty_{-1}, \Sigma j]_{s+t}
$$
all the groups $E^2_{s, t}$ with $s+t=0$ are zero.
\qed
\enddemo

\remark{Remark 3.7}
Let $G/O$ be the homotopy fiber of the map of spaces $BO \to BG$, and
let $\Cok J = \Omega^\infty c$ be the cokernel-of-J space.  There is a
(Sullivan) fiber sequence $\Cok J \to G/O \to BSO$ [\Cite{MM}, \S5C].
Waldhausen [\Cite{Wa3}, 3.4] constructed a space level map $hw \: G/O
\to \Omega \Omega^\infty Wh(*)$, using manifold models for $A(*)$.
Hence there is a geometrically defined composite map $\Cok J \to G/O
\to \Omega \Omega^\infty Wh(*)$.  Presumably this is homotopic to the
infinite loop map $\Omega^\infty \Sigma^{-1} g$.
\endremark

\subhead A cofiber sequence \endsubhead
We can analyze a variant of the lower cofiber sequence in~3.6 by passing
to connective covers.  There is a map of homotopy Cartesian squares from
$$
\xymatrix{
\Sigma c \ar[r]^-{\Sigma f} \ar[d]^{*} &
S^1 \ar[d]^{\In_2 \circ \Sigma e} \\
j \ar[r]^-{\In_1} &
j \vee \Sigma j
}
\qquad\text{to}\qquad
\xymatrix{
Wh(*) \ar[r]^-{\widetilde{\trc}} \ar[d]^{L\iota} &
P_0 \Sigma \C P^\infty_{-1} \ar[d]^{L\iota} \\
K(\Z) \ar[r]^-{\trc_{\Z}} &
P_0 TC(\Z; p)
}
$$
induced by $g$, $\Sigma b_0$, $\In_1$ and $\In_1 \vee \In_2$ in the upper
left, upper right, lower left and lower right corners, respectively.
In the lower rows we are using the splittings $K(\Z) \simeq j \vee
\Sigma^5 ko$ and $P_0 TC(\Z; p) \simeq K(\Z_p) \simeq j \vee \Sigma
j \vee \Sigma^3 ku$ derived from~3.1.  Let $\tau \: Wh(*)/\Sigma c
\to P_0 \Sigma\overline{\C P}^\infty_{-1}$, $\ell \: Wh(*)/\Sigma c
\to \Sigma^5 ko$ and $\ell \: P_0 \Sigma\overline{\C P}^\infty_{-1}
\to \Sigma^3 ku$ be the cofiber maps induced by $\widetilde{\trc} \:
Wh(*) \to P_0 \Sigma\C P^\infty_{-1}$, $L\iota \: Wh(*) \to K(\Z)$ and
$L\iota \: P_0 \Sigma \C P^\infty_{-1} \to P_0 TC(\Z; p)$, respectively.

\proclaim{Theorem 3.8}
Assume~3.2.
There is a diagram of horizontal and vertical cofiber sequences:
$$
\xymatrix{
&
{}\hofib(\ell) \ar[r]^{=} \ar[d] &
{}\hofib(\ell) \ar[d] \\
\Sigma^2 ko \ar[r] \ar[d]^{=} &
Wh(*)/\Sigma c \ar[r]^\tau \ar[d]^\ell &
P_0 \Sigma\overline{\C P}^\infty_{-1} \ar[r] \ar[d]^\ell &
\Sigma^3 ko \ar[d]^{=} \\
\Sigma^2 ko \ar[r]^{*} &
\Sigma^5 ko \ar[r]^c &
\Sigma^3 ku \ar[r]^{r\beta^{-1}} &
\Sigma^3 ko \,.
}
$$
The map $\tau \: Wh(*)/\Sigma c \to P_0 \Sigma\overline{\C P}^\infty_{-1}$
induces a split injection on homotopy groups in all degrees, and each
map $\ell$ is $(2p-3)$-connected.  Thus
$$
\pi_*(\tau) \: \tors \pi_*(Wh(*)/\Sigma c) \cong
\tors \pi_*(\Sigma\overline{\C P}^\infty_{-1}) \,.
$$
\endproclaim

Here $\tors A$ denotes the torsion subgroup of an abelian group $A$.

\demo{Proof}
It follows from~3.1 and localization in algebraic $K$-theory that the
map $\Sigma^5 ko \to \Sigma^3 ku$ induced by $\trc_{\Z} \: K(\Z) \to
P_0 TC(\Z; p) \simeq K(\Z_p)$ is the lift of $\Sigma c \: \Sigma ko
\to \Sigma ku$ to the 1-connected covers.  This identifies the central
homotopy Cartesian square in the diagram.

By comparing the vertical homotopy fibers in the last three homotopy
Cartesian squares we obtain a cofiber sequence $c \vee \Sigma c \to
\hofib(L) \to \hofib(\ell)$, as in [\Cite{KR}, 3.6].  Hence each map $\ell$
is $(2p-3)$-connected because $L$ is.  There is a $(4p-3)$-connected space
level map from $SU$ to $\Omega^\infty \Sigma \overline{\C P}^\infty_{-1}$,
as in [\Cite{KR}, (17)].
$$
B\phi \: SU @>>> \Omega^\infty \Sigma \overline{\C P}^\infty_{-1}
@>\Omega^\infty\ell>> SU \,.
$$
Its composite with $\Omega^\infty\ell$ to $\Omega^\infty \Sigma^3 ku =
SU$ loops to an H-map $\phi \: BU \to BU$.  Any such H-map is a series
of Adams operations $\psi^k$, as in [\Cite{MST}, 2.3], so $\pi_*(\phi;
\Z/p)$ only depends on $* \bmod q$ in positive degrees.  Since $\ell$
is $(2p-3)$-connected it follows that $\phi$ is $(2p-4)$-connected,
so $\pi_*(\phi; \Z/p)$ is an isomorphism for $0 < * < q$, and so
$\pi_*(\phi)$ is an isomorphism for all $* \not\equiv 0 \bmod q$.
Hence $\pi_*(\ell)$ is (split) surjective whenever $* \not\equiv 1 \bmod
q$, cf.~[\Cite{KR}, 6.3(i)].

Finally $r\beta^{-1}$ is split surjective as a spectrum map, and
$\pi_*(\Sigma^3 ko)$ is zero for $* \equiv 1 \bmod q$, so $r\beta^{-1}\ell
\: P_0 \Sigma\overline{\C P}^\infty_{-1} \to \Sigma^3 ko$ induces a
split surjection on homotopy in all degrees.
\qed
\enddemo

\remark{Remark 3.9}
We still do not know the behavior of $\ell \: Wh(*)/\Sigma c \to
\Sigma^5 ko$ in degrees $* \equiv 1 \bmod q$.  It induces the same
homomorphism on homotopy as $\ell \: P_0 \Sigma\overline{\C P}^\infty_{-1} \to \Sigma^3 ku$, since
$\pi_*(\tau)$ and $\pi_*(c)$ are isomorphisms in these degrees.
\endremark

\remark{Remark 3.10}
By a result of Madsen and Schlichtkrull [\Cite{MS}, 1.3] there
is a splitting of implicitly $p$-completed spaces $\Omega^\infty
(\Sigma\overline{\C P}^\infty_{-1}) \simeq Y \times SU$, where $\pi_*(Y)
\cong \tors \pi_*(\Sigma\overline{\C P}^\infty_{-1})$ is finite in
each degree.  The map
$$
Y \times SU \simeq \Omega^\infty (\Sigma\overline{\C P}^\infty_{-1})
@>\Omega^\infty (r\beta^{-1}\ell)>> \Omega^\infty (\Sigma^3 ko) \simeq
Sp \simeq SO
$$
induces a split surjection on homotopy
groups in all degrees, so the composite map $SU @>\In_2>> Y \times SU \to SO$
has homotopy fiber $BBO$, by real Bott periodicity.  Hence there is a fiber
sequence
$$
BBO \to \Omega^\infty (Wh(*)/\Sigma c) \to Y
$$
and split short exact sequences
$$
0 \to \pi_*(BBO) \to \pi_*(Wh(*)/\Sigma c) \to \pi_*(Y) \to 0
$$
in each degree.
\endremark

\subhead The suspended quaternionic projective spectrum \endsubhead
After $p$-adic completion $\C P^\infty_{-1}$ splits as a wedge sum of
$(p-1)$ eigenspectra $\C P^\infty_{-1}[a]$ for $-1 \le a \le p-3$,
much like the $p$-complete (or $p$-local) Adams splitting of $ku$
from \cite{Ad}, and the $p$-complete splitting of $\Sigma^\infty(\C
P^\infty_+)$ from [\Cite{MT}, \S4.1].  Here $H^*(\C P^\infty_{-1}[a])
\cong \Z_p\{y^k \mid k \ge -1, k \equiv a \bmod p-1\}$, and similarly
with mod~$p$ coefficients.

Let $\H P^\infty$ be the infinite quaternionic projective spectrum.
The ``quaternionification'' map $q \: \C P^\infty_{-1} \to \H P^\infty_+
\simeq S \vee \H P^\infty$ admits a (stable $p$-adic) section $s
\: \H P^\infty_+ \to \C P^\infty_{-1}$.  (It can be obtained by
Thomifying the Becker--Gottlieb transfer map $\Sigma^\infty(BS^3_+)
\to \Sigma^\infty(BS^1_+)$ associated to the sphere bundle $S^2 \to
BS^1 \to BS^3$, with respect to minus the tautological quaternionic line
bundle over $BS^3 = \H P^\infty$, and collapsing the bottom $(-4)$-cell.
It is a section because the Euler characteristic $\chi(S^2) = 2$ is a
unit mod~$p$.)  This section $s$ identifies $S \vee \H P^\infty$ with
the wedge sum of the even summands $\C P^\infty_{-1}[a]$ for $a = 2i$
with $0 \le i \le (p-3)/2$.

Splitting off $S$, suspending once and passing to connected covers,
we obtain maps $s' \: \Sigma \H P^\infty \to P_0 \Sigma\overline{\C
P}^\infty_{-1}$ and $q' \: P_0 \Sigma\overline{\C P}^\infty_{-1} \to
\Sigma \H P^\infty$ whose composite is a $p$-adic equivalence.

\proclaim{Proposition 3.11}
The map $s' \: \Sigma \H P^\infty \to P_0 \Sigma\overline{\C
P}^\infty_{-1}$ admits a lift
$$
\tilde s \: \Sigma \H P^\infty \to Wh(*)/\Sigma c
$$
over $\tau$, which is unique up to homotopy, and whose composite with
$$
q' \circ \tau \: Wh(*)/\Sigma c \to \Sigma \H P^\infty
$$
is a $p$-adic equivalence.
\endproclaim

\demo{Proof}
The composite map $r\beta^{-1}\ell \circ s' \: \Sigma \H P^\infty
\to \Sigma^3 ko$ lies in a zero group, by the Atiyah--Hirzebruch spectral
sequence
$$
E^2_{s, t} = H^{-s}(\Sigma \H P^\infty; \pi_t \Sigma^3 ko)
\Longrightarrow [\Sigma \H P^\infty, \Sigma^3 ko]_{s+t} \,.
$$
Hence $s'$ admits a lift $\tilde s$, as claimed.
In fact the lift is unique up to homotopy, since also $[\Sigma \H
P^\infty, \Sigma^3 ko]_1 = 0$.
\qed
\enddemo

\subhead A second cofiber sequence \endsubhead
We define $Wh(*)/(\Sigma c, \Sigma \H P^\infty) \simeq \hofib(q'\tau)$
as the homotopy cofiber of $\tilde s$, and write
$$
{P_0 \Sigma\overline{\C P}^\infty_{-1} \over \Sigma \H P^\infty}
\simeq P_0 \Sigma \C P^\infty_{-1}[-1] \vee \bigvee_{i=1}^{(p-3)/2}
\Sigma \C P^\infty_{-1}[2i-1]
\tag 3.12
$$
for the suspended homotopy cofiber of $s'$. Then:

\proclaim{Theorem 3.13}
Assume~3.2.
There is a splitting
$$
Wh(*) \simeq \Sigma c \vee \Sigma \H P^\infty
\vee {Wh(*) \over (\Sigma c, \Sigma \H P^\infty)}
$$
and a cofiber sequence
$$
\Sigma^2 ko @>>>
{Wh(*) \over (\Sigma c, \Sigma \H P^\infty)} @>\tau>>
{P_0 \Sigma\overline{\C P}^\infty_{-1} \over \Sigma \H P^\infty} @>\delta>>
\Sigma^3 ko \,.
$$
The map $\tau$ induces a split injection on homotopy groups in all
degrees, and the map $\delta$ induces an injection on mod~$p$ cohomology
in degrees $\le 2p-3$.  Thus
$$
\pi_*(Wh(*)) \cong \pi_*(\Sigma c) \oplus \pi_*(\Sigma \H P^\infty)
\oplus \tors \pi_*\biggl({\Sigma\overline{\C P}^\infty_{-1} \over
\Sigma \H P^\infty}\biggr) \,.
$$
\endproclaim

\demo{Proof}
The cofiber sequence arises by splitting off $\Sigma \H P^\infty$
from the middle horizontal cofiber sequence in~3.8.  The assertion
about $\tau$ follows by retraction from the corresponding
statement in~3.8.  The map $\delta$ is the composite of the maps
$$
{P_0 \Sigma\overline{\C P}^\infty_{-1} \over \Sigma \H P^\infty} @>\In>>
P_0 \Sigma\overline{\C P}^\infty_{-1} @>\ell>>
\Sigma^3 ku @>r\beta^{-1}>> \Sigma^3 ko \,.
$$
On mod~$p$ cohomology $(r\beta^{-1})^*$ is split injective
and $\ell^*$ is injective in degrees $\le 2p-3$ by~3.8.  The 
kernel of $\In^*$ is $\Sigma H^*(\H P^\infty; \F_p)$, which is
concentrated in degrees $\equiv 1 \bmod 4$.  But in degrees $\le
2p-3$ all of $H^*(\Sigma^3 ko; \F_p)$ is in degrees $\equiv 3 \bmod
4$, so also the composite $\delta^*$ is injective in this range
of degrees.
\qed
\enddemo

\remark{Remark 3.14}
Note that the upper cofiber sequence in~3.4 maps as in~3.6 to the middle
horizontal cofiber sequence in~3.8, which in turn maps to the cofiber
sequence in~3.13.  In~5.4 we will see that $\delta$ is $(4p-2)$-connected.
\endremark

\head Homotopical analysis \endhead

\subhead Homotopy of the fiber of the restricted $S^1$-transfer
map \endsubhead
To make the $p$-primary homotopy groups of $Wh(*)$ explicit we refer
to~3.8 and compute the $p$-torsion in the homotopy of $\overline{\C
P}^\infty_{-1}$ in an initial range of degrees.  This is related to $\C
P^\infty$ by the cofiber sequence
$$
\Sigma\overline{\C P}^\infty_{-1} @>>> \Sigma \C P^\infty @>t>> S
\tag 4.1
$$
extracted from~(3.5).  We also use the cofiber sequence
$$
c \wedge \C P^\infty @>f\wedge1>> \C P^\infty @>e\wedge1>> j \wedge
\C P^\infty
$$
obtained by smashing (2.5) with $\C P^\infty$.  There are
Atiyah-Hirzebruch spectral sequences:
$$
\align
E^2_{s, t} &= H_s(\C P^\infty; \pi_t(j)) \Longrightarrow j_{s+t}(\C
P^\infty)
\tag 4.2 \\
E^2_{s, t} &= H_s(\C P^\infty; \pi_t(S)) \Longrightarrow \pi_{s+t}(\C
P^\infty)
\tag 4.3 \\
E^2_{s, t} &= H_s(\overline{\C P}^\infty_{-1}; \pi_t(S))
\Longrightarrow \pi_{s+t}(\overline{\C P}^\infty_{-1}) \,.
\tag 4.4 \\
\endalign
$$
We will now account for the abutment of~(4.2) in all degrees, and
for~(4.3) and~(4.4) in total degrees $* < |\beta_2 b_1| = (2p+1)q$ and $*
< |\beta_2 b_{-1}| = (2p+1)q-4$, respectively.

Let $v_p(n)$ be the $p$-adic valuation of a natural number $n$.
In degrees $* < |\beta_2| = (2p+1)q-2$ the $p$-torsion in $\pi_*(S)
= \pi^S_*$ is generated by the image-of-J classes $\bar\alpha_i \in
\pi^S_{qi-1}$ of order $p^{1+v_p(i)}$ for $i\ge1$, and the cokernel-of-J
classes [\Cite{Ra}, 1.1.14]
$$
\beta_1 \in \pi^S_{pq-2}, \quad \alpha_1 \beta_1 \in \pi^S_{(p+1)q-3},
\quad \beta_1^2 \in \pi^S_{2pq-4} \quad\text{and}\quad \alpha_1 \beta_1^2
\in \pi^S_{(2p+1)q-5} \,,
$$
each of order~$p$.

\proclaim{Theorem 4.5}
Above the horizontal axis and in total degrees $* < |\beta_2|-2$,
the Atiyah--Hirzebruch $E^\infty_{s, t}$-term for $\pi_* \overline{\C
P}^\infty_{-1}$ agrees with that for $j_*(\C P^\infty)$, {\bf plus} the
$\Z/p$-module generated by $\beta_1 b_m$, $\alpha_1 \beta_1 b_{mp}$,
$\beta_1^2 b_m$ \rom(and $\alpha_1 \beta_1^2 b_{mp}$, which is in a higher
total degree\rom) for $1 \le m \le p-3$, {\bf minus} the $\Z/p$-module
generated by $\alpha_1 b_{mp}$ for $m \ge p-2$.
\endproclaim

We give the proof in a couple of steps.

\subhead Connective J-theory of complex projective space \endsubhead
On the horizontal axis the $E^2$-terms of~(4.2) and~(4.3) have the
form $E^2_{*, 0} = H_*(\C P^\infty) = \Z\{b_n \mid n\ge1\}$, which has
the structure of a divided power algebra on $b_1$.  By Toda \cite{To}
or Mosher [\Cite{Mo}, 2.1], the corresponding part of the $E^\infty$-term
of~(4.3) consists of the polynomial algebra on $b_1$, i.e.,
$$
E^\infty_{2n, 0} = \Z\{n! \, b_n\} \subseteq E^2_{2n, 0} = \Z\{b_n\}
\tag 4.6
$$
for all $n\ge1$.  Hence the order of the images of the differentials
$d^r_{2n, 0}$ landing in total degree~$2n-1$ all multiply to $n!$.

It is known by [\Cite{Mo}, 4.7(a)] that these differentials from the
horizontal axis land in the image-of-J, i.e., have the form $\theta
b_k$ with $\theta$ a multiple of some $\bar\alpha_i$.  Hence~(4.6) also
gives the $E^\infty$-term of~(4.2) on the horizontal axis.  Since the
Atiyah--Hirzebruch spectral sequence for $j_*(\C P^\infty)$ only has
classes in (even, odd) bidegrees above the horizontal axis, there can
be no further differentials in~(4.2).  In even total degrees it follows
that $j_{2n}(\C P^\infty) \cong \Z\{n! \, b_n\}$ for $n\ge1$.

In odd total degrees, the $E^2$-term of~(4.2) contains the classes
$p^e \bar\alpha_i b_k$ in bidegree $(s, t) = (2k, qi-1)$, for $0
\le e \le v_p(i)$.  It follows that the $p$-valuation of the order of
the groups $E^2_{s, t}$ in total degree $s+t=2n-1$ equals $\sum_{e\ge0}
\left[ (n-1) / p^e(p-1) \right]$, so the $p$-valuation of the order of
the finite group $j_{2n-1}(\C P^\infty)$ is
$$
\sum_{e\ge0} \left[ n-1 \over p^e(p-1) \right]
-\sum_{e\ge0} \left[ n \over p^e \cdot p \right] \,.
$$
Here the second sum equals $v_p(n!)$.  Compare [\Cite{KR}, 4.3]
due to Knapp.  For $n \le p^2(p-1)$ the terms with $e\ge2$ vanish.

\comment
for $e=0$ we have
$$
\left[ n-1 \over p-1 \right] - \left[ n \over p \right]
= \left[ n-1 \over p(p-1) \right] + \cases
1 & \text{for $m(p-1) < n < mp$, $1 < m < p$,} \\
0 & \text{otherwise,}
\endcases
$$
and for $e=1$ we have
$$
\left[ n-1 \over p(p-1) \right] - \left[ n \over p^2 \right]
= \cases
1 & \text{for $p^2-p+1 \le n \le p^2-1$,} \\
0 & \text{for $1 \le n \le p^2-p$ or $p^2 \le n \le 2p^2-2p$.} \\
\endcases
$$
\endcomment

\subhead Stable homotopy of complex projective space \endsubhead
We now return to~(4.3) where the $E^2$-term contains additional classes
from $H_*(\C P^\infty; \pi_*(c))$.  The primary operation $P^1$ detects
$\alpha_1$, and $P^1(y^k) = k y^{k+p-1}$ in mod~$p$ cohomology, so there
are differentials $d^q(\theta b_{k+p-1}) = k\alpha_1\theta b_k$ for all
$\theta \in \pi_*(S)$.  In the case $\theta=1$ these differentials were
already accounted for by the differentials leading to~(4.6), but for $t <
|\beta_2|$ there are also differentials
$$
d^q(\beta_1 b_{k+p-1}) = \alpha_1\beta_1 b_k
\qquad\text{and}\qquad
d^q(\beta_1^2 b_{k+p-1}) = \alpha_1\beta_1^2 b_k
$$
up to unit multiples, for $k \not\equiv 0 \bmod p$, $k\ge1$.  This
leaves the classes $\alpha_1 b_{mp}$ (already in $j_*(\C P^\infty)$),
$\alpha_1\beta_1 b_{mp}$ and $\alpha_1\beta_1^2 b_{mp}$ for $m\ge1$ in
odd total degrees, and the classes $\beta_1 b_1, \dots, \beta_1 b_{p-2}$,
$\beta_1 b_{mp-1}$ for $m\ge1$, $\beta_1^2 b_1, \dots, \beta_1^2 b_{p-2}$
and $\beta_1^2 b_{mp-1}$ for $m\ge1$ in even total degrees.

The (well-known) $p$-fold Toda bracket $\beta_1 = \langle \alpha_1,
\dots, \alpha_1 \rangle$ implies differentials
$$
d^{(p-1)q}(\theta \alpha_1 b_{k+(p-1)^2}) = \theta \beta_1 b_k
$$
when $k+(p-1)^2 = mp$, up to unit multiples.  So the classes $\alpha_1
b_{mp}$ (from $j_*(\C P^\infty)$) and $\alpha_1\beta_1 b_{mp}$ for $m
\ge p-1$ support $d^{(p-1)q}$-differentials, which kill the classes
$\beta_1 b_{mp-1}$ and $\beta_1^2 b_{mp-1}$ for $m\ge1$.  For bidegree
reasons this accounts for all differentials in~(4.3) in total degrees $*
< |\beta_2 b_1|$.

To pass from $\C P^\infty$ to $\overline{\C P}^\infty_{-1}$ we must take
into account the differentials in~(4.4) that cross the vertical axis,
which amounts to the restricted $S^1$-transfer map $t$ as in~(4.1).  The
image-of-J in its target $\pi_*(S)$ is hit by classes on the horizontal
axis of~(4.3), by [\Cite{Mu}, 4.3] or Crabb and Knapp, cf.~[\Cite{KR}, 5.8].
The cokernel-of-J classes are hit by the differentials
$$
\alignat 2
d^q(\beta_1 b_{p-2}) &= \alpha_1\beta_1 b_{-1} \,,
&d^{(p-1)q}(\alpha_1 b_{(p-2)p}) &= \beta_1 b_{-1} \,, \\
d^q(\beta_1^2 b_{p-2}) &= \alpha_1\beta_1^2 b_{-1} \,, \qquad
&d^{(p-1)q}(\alpha_1\beta_1 b_{(p-2)p}) &= \beta_1^2 b_{-1}
\endalignat
$$
in~(4.4).
Looking over the bookkeeping concludes the proof of Theorem~4.5.

\subhead Torsion homotopy of the smooth Whitehead spectrum \endsubhead
\proclaim{Theorem 4.7}
{\rm(a)}\qua
Assume~3.2.
The torsion homotopy of $Wh(*)$ decomposes as
$$
\tors \pi_*(Wh(*)) \cong \pi_*(\Sigma c) \oplus \tors \pi_*(\Sigma
\overline{\C P}^\infty_{-1})
$$
in all degrees.

{\rm(b)}\qua
In degrees $* < |\beta_2|+1 = (2p+1)q-1$
$$
\pi_*(\Sigma c) \cong \Z/p\{ \Sigma\beta_1, \Sigma\alpha_1\beta_1,
\Sigma\beta_1^2, \Sigma\alpha_1\beta_1^2 \}
$$
with generators in degrees $pq-1$, $(p+1)q-2$, $2pq-3$ and $(2p+1)q-4$,
respectively.

{\rm(c)}\qua
In even degrees $* < |\beta_2|-1 = (2p+1)q-3$ the $p$-valuation of the
order of $\tors \pi_{2n}(\Sigma \overline{\C P}^\infty_{-1})$ equals
$$
\left(\left[ n-1 \over p-1 \right] + \left[ n-1 \over p(p-1) \right] \right)
-
\left(\left[ n \over p \right] + \left[ n \over p^2 \right]\right)
\,,
$$
{\bf plus}~$1$ when $n = p^2-2+mp$ for $1 \le m \le p-3$,
{\bf minus}~$1$ when $n = p-1+mp$ for $m \ge p-2$.

{\rm(d)}\qua
In odd degrees $* < |\beta_2|-1 = (2p+1)q-3$ the $p$-valuation of the
order of $\tors \pi_{2n+1}(\Sigma \overline{\C P}^\infty_{-1})$ equals
$1$ when $n = p^2-p-1+m$ or $n = 2p^2-2p-2+m$ for $1 \le m \le
p-3$, and is~$0$ otherwise.
\endproclaim

\example{Example 4.8}
(a)\qua
When $p=3$, the 3-torsion in $\pi_* Wh(*)$ has order~3 in degrees~11, 16,
18, 20, 21 and~22, order~$3^2$ in degree~24, order~$3^3$ in degree~14,
and is trivial in the remaining degrees $*<25$.

(b)\qua
When $p=5$, the 5-torsion in $\pi_* Wh(*)$ has order~5 in degrees~18,
26, 28, 34, 36, 39, 41, 43, 48, 50, 52, 54, 58, 60, 62, 64, 68, 70, 72,
77, 78, 79, 80 and~81, order~$5^2$ in degrees~42, 44, 56, 74 and~76,
order~$5^3$ in degrees~46, 66 and~82, order~$5^4$ in degree~84, and is
trivial in the remaining degrees $*<85$.
\endexample

In roughly half this range we can give the following simpler statement.

\proclaim{Corollary 4.9}
{\rm(a)}\qua
For $p\ge5$, the low-degree $p$-torsion in $\pi_* Wh(*)$ is $\Z/p$
in degrees $* = 2n$ for $m(p-1) < n < mp$ and $1 < m < p$, except in
degree~$2p^2-2p-2$ \rom(corresponding to $n = mp-1$ and $m=p-1$\rom).
The next $p$-torsion is $\Z/p\{\Sigma\beta_1\}$ in degree $2p^2-2p-1$,
and a group of order~$p^2$ in degree $2p^2-2p+2$.

{\rm(b)}\qua
For $p=3$ the bottom $3$-torsion in $\pi_* Wh(*)$ is
$\Z/3\{\Sigma\beta_1\}$ in degree~$11$, followed by
$\Z/3\{\Sigma\alpha_1\beta_1\} \oplus \Z/9$ in degree~$14$.
\endproclaim

The asserted group structure of $\pi_{14} Wh(*)_{(3)}$ can be obtained
from 5.5(a) below and the mod~3 Adams spectral sequence.

\remark{Remark 4.10}
Klein and the author showed in [\Cite{KR}, 1.3(iii)] that for any odd
prime~$p$, regular or irregular, below degree $2p^2-2p-2$ there are
direct summands $\Z/p$ in $\pi_{2n} Wh(*)$ for $m(p-1) < n < mp$ and $1 <
m < p$.  The calculations above show that under the added hypothesis~3.2,
these classes constitute all of the $p$-torsion in $\pi_* Wh(*)$, in
this range of degrees.
\endremark

\head Cohomological analysis \endhead

We can determine the mod~$p$ cohomology of $Wh(*)$ as a module over the
Steenrod algebra~$A$, up to an extension, in all degrees.  To do this,
we apply cohomology to the splitting and cofiber sequence in~3.13.

\subhead Some cohomology modules \endsubhead
Let us briefly write $H^*(X) = H^*(X; \F_p)$ for the mod~$p$ cohomology of
a spectrum $X$, where $p$ is an odd prime.  It is naturally a left module
over the mod~$p$ Steenrod algebra $A$ \cite{St}.  Let $A_n$ be the
subalgebra of $A$ generated by the Bockstein operation $\beta$ and the
Steenrod powers $P^1, \dots, P^{p^{n-1}}$ and let $E_n$ be the exterior
subalgebra generated by the Milnor primitives $\beta, Q_1, \dots, Q_n$,
where $Q_0 = \beta$ and $Q_{n+1} = [P^{p^n}, Q_n]$.  For an augmented
subalgebra $B \subset A$ we write $I(B) = \ker(\epsilon \: B \to \F_p)$
for the augmentation ideal, and let $A//B = A \otimes_B \F_p = A/A\cdot I(B)$.

\proclaim{Proposition 5.1}
{\rm(a)}\qua $H^*(H\Z) \cong A//E_0 = A/A(\beta)$ and $H^*(\ell) \cong A//E_1
= A/A(\beta, Q_1)$.

{\rm(b)}\qua The cofiber sequence $\Sigma^{q-1} \ell \to j \to \ell$ induces a
nontrivial extension
$$
0 \to A//A_1 \to H^*(j) \to \Sigma^{pq-1} A//A_1 \to 0
$$
of $A$-modules.  As an $A$-module $H^*(j)$ is generated by two classes $1$
and $b$ in degree~$0$ and $pq-1$, respectively, with $\beta(b) = P^p(1)$.

{\rm(c)}\qua The cofiber sequence $S @>e>> j \to \Sigma c$ induces an
identification $H^*(\Sigma c) \cong \ker(e^* \: H^*(j) \to \F_p)$.
There is a nontrivial extension
$$
0 \to I(A)/A(\beta, P^1) \to H^*(\Sigma c) \to \Sigma^{pq-1} A//A_1 \to 0
$$
of $A$-modules.
\endproclaim

\demo{Proof}
For~(a), see [\Cite{AP}, 2.1].
For~(c), clearly the given cofiber sequence identifies $H^*(\Sigma c)$
with the positive degree part of $H^*(j)$.
The long exact sequence in cohomology associated
to the cofiber sequence given in~(b) is:
$$
\Sigma^q A//E_1 @>(\psi^r-1)^*>> A//E_1 @>>> H^*(j) @>>>
\Sigma^{q-1} A//E_1 @>(\psi^r-1)^*>> \Sigma^{-1} A//E_1 \,.
$$
The map $e \: S \to j$ is $(pq-2)$-connected [\Cite{Ra}, 1.1.14],
so $e^* \: H^*(j) \to H^*(S) = \F_p$ is an isomorphism for $* \le pq-2$.
Thus $P^1 \in A//E_1$ is in the image of $(\psi^r-1)^*$, and so
$(\psi^r-1)^*$ is induced up over $A_1 \subset A$ by
$$
\Sigma^q A_1//E_1 @>P^1>> A_1//E_1 \,,
$$
which has kernel $\Sigma^{pq}\F_p$ generated by $\Sigma^q P^{p-1}$ and
cokernel $\F_p$ generated by $1$.  Hence there is an extension $A//A_1
\to H^*(j) \to \Sigma^{pq-1} A//A_1$.  Note that the bottom classes
in $A//A_1$ are $1$ and $P^p$ in degrees~0 and $pq$, respectively.
Let $b \in H^{pq-1}(j)$ be the class mapped to $\Sigma^{pq-1}(1)$
in $\Sigma^{pq-1} A//A_1$.  By the Hurewicz theorem for $\Sigma c$
it is dual to the Hurewicz image of the bottom class $\Sigma \beta_1
\in \pi_{pq-1}(\Sigma c)$.  Since $\beta_1 \in \pi_{pq-2}(c) \subset
\pi_{pq-2}(S)$ has order~$p$ there is a nontrivial Bockstein $\beta(b)$
in $H^*(\Sigma c)$, and thus also in $H^*(j)$.  The only possible value
in degree $pq$ is $P^p(1)$.  Part~(c) now follows easily from~(b).
\qed
\enddemo

\proclaim{Proposition 5.2}
{\rm(a)}\qua 
$H^*(\Sigma \H P^\infty) \cong \F_p\{\Sigma y^k
\mid \text{$k\ge2$ even}\}$.

{\rm(b)}\qua
$H^*(\Sigma \C P^\infty_{-1}[-1]) \cong \Sigma^{-1} A/C$.  Here $C
\subset A$ is the annihilator ideal of $\Sigma y^{-1}$, which is spanned
over $\F_p$ by all admissible monomials in $A$ except $1$ and the $P^i$
for $i\ge1$.

{\rm(c)}\qua
The cofiber sequence $P_0 \Sigma \C P^\infty_{-1}[-1] \to \Sigma \C
P^\infty_{-1}[-1] \to \Sigma^{-1} H\Z$ induces an identification $H^*(P_0
\Sigma \C P^\infty_{-1}[-1]) \cong \Sigma^{-2} C/A(\beta)$.

{\rm(d)}\qua
For $1 \le i \le (p-3)/2$ there are isomorphisms $H^*(\Sigma \C
P^\infty_{-1}[2i-1]) \cong \F_p\{\Sigma y^k \mid k = 2i-1 + m(p-1),
m \ge 0\}$.
\endproclaim

\demo{Proof}
Any admissible monomial $P^I$ with $I = (i_1, \dots, i_n)$ and
$n\ge2$ acts trivially on $\Sigma y^{-1}$ because $z = P^{i_n}(\Sigma
y^{-1})$ is in the image from $H^*(\Sigma \C P^\infty)$, which is an
unstable $A$-module, and then $P^{i_{n-1}}(z) = 0$ by instability.
\qed
\enddemo

\subhead Cohomology of the smooth Whitehead spectrum \endsubhead
\proclaim{Proposition 5.3}
The $A$-module homomorphism
$$
\delta^* \: H^*(\Sigma^3 ko) \to H^*(P_0 \Sigma\overline{\C
P}^\infty_{-1}/\Sigma \H P^\infty)
$$
splits as the direct sum of the injection
$$
\Sigma^{q-1} A//E_1 @>>> \Sigma^{-2} C/A(\beta)
$$
taking $\Sigma^{q-1}(1)$ to $\Sigma^{-2} Q_1$,
and the homomorphisms
$$
\align
\delta_i^* \: \Sigma^{4i-1} A//E_1 &@>>> H^*(\Sigma \C P^\infty_{-1}[2i-1]) \\
&\cong \F_p\{\Sigma y^k \mid k = 2i-1 + m(p-1), m \ge 0\}
\endalign
$$
taking $\Sigma^{4i-1}(1)$ to $\Sigma y^{2i-1}$ for $1 \le i \le (p-3)/2$.
\endproclaim

\demo{Proof}
By~(2.4) and~5.1(a) the source of $\delta^*$ splits as the direct sum of
the cyclic $A$-modules $\Sigma^{4i-1} A//E_1$ for $1 \le i \le (p-1)/2$.
Here $4i-1 = q-1$ for $i = (p-1)/2$.  Hence $\delta^*$ is determined as an
$A$-module homomorphism by its value on the generators $\Sigma^{4i-1}(1)$.
These are all in degrees $\le q-1=2p-3$, and $\delta^*$ is injective in
this range by~3.13.  By~(3.12), 5.2(c) and~(d) the target of $\delta^*$
splits as the direct sum of $\F_p\{\Sigma y^k \mid k \equiv 2i-1 +
m(p-1), m\ge0\}$ for $1 \le i \le (p-3)/2$ and $\Sigma^{-2} C/A(\beta)$.
The bottom class of the latter is $\Sigma^{-2} Q_1$, in degree $q-1$.
Hence the target of $\delta^*$ has rank~1 in each degree $4i-1$ for $1
\le i \le (p-1)/2$, and so (up to a unit which we suppress) $\delta^*$
maps $\Sigma^{4i-1}(1)$ to $\Sigma y^{2i-1}$ for $1 \le i \le (p-3)/2$
and $\Sigma^{q-1}(1)$ to $\Sigma^{-2} Q_1$.

The homomorphism $\Sigma^{q-1} A//E_1 \to \Sigma^{-2} C/A(\beta)$ is
injective, as its continuation into $\Sigma^{-2} A//E_0$ is induced up
over $E_1 \subset A$ from the injection $\Sigma^{q-1} \F_p \to \Sigma^{-2}
E_1//E_0$ taking $\Sigma^{q-1}(1)$ to $\Sigma^{-2} Q_1$.
\qed
\enddemo

\proclaim{Theorem 5.4}
Assume~3.2.
There is a splitting
$$
H^*(Wh(*)) \cong H^*(\Sigma c) \oplus H^*(\Sigma \H P^\infty) \oplus
H^*\biggl({Wh(*) \over (\Sigma c, \Sigma \H P^\infty)}\biggr)
$$
and an extension of $A$-modules
$$
0 \to \cok(\delta^*) @>>> H^*\biggl({Wh(*) \over (\Sigma c, \Sigma \H
P^\infty)}\biggr) @>>> \Sigma^{-1} \ker(\delta^*) \to 0
$$
where
$$
\cok(\delta^*) \cong \Sigma^{-2} C/A(\beta, Q_1) \oplus
\bigoplus_{i=1}^{(p-3)/2} H^*(\Sigma \C P^\infty_{-1}[a])/A(\Sigma
y^a)
$$
and
$$
\Sigma^{-1} \ker(\delta^*) \cong \bigoplus_{i=1}^{(p-3)/2} \Sigma^{2a}
C_a/A(\beta, Q_1) \,.
$$
In both sums we briefly write $a = 2i-1$, so $a$ is odd with $1 \le a
\le p-4$.  Here $H^*(\Sigma \C P^\infty_{-1}[a]) = \F_p\{\Sigma y^k
\mid k \equiv a \bmod p-1, k \ge a\}$, $A(\Sigma y^a) \subset H^*(\Sigma
\C P^\infty_{-1}[a])$ is the submodule generated by $\Sigma y^a$, and
$C_a \subset A$ is the annihilator ideal of $\Sigma y^a \in H^*(\Sigma
\C P^\infty_{-1}[a])$.
\endproclaim

\demo{Proof}
The splitting and extension follow by applying cohomology to~3.13.
The cohomologies of $\Sigma c$ and $\Sigma \H P^\infty$ are given
in~5.1(c) and~5.2(a), respectively.  The descriptions of $\ker(\delta^*)$
and $\cok(\delta^*)$ are immediate from~5.3.
\qed
\enddemo

\example{Example 5.5}
(a)\qua  When $p=3$ there is a splitting
$$
H^*(Wh(*)) \cong H^*(\Sigma c) \oplus H^*(\Sigma \H P^\infty)
\oplus \Sigma^{-2} C/A(\beta, Q_1) \,.
$$

(b)\qua When $p=5$ there is an extension
$$
\multline
0 \to \Sigma^{-2} C/A(\beta, Q_1) \oplus H^*(\Sigma \C
P^\infty_{-1}[1])/A(\Sigma y) \\
@>>> H^*\biggl({Wh(*) \over (\Sigma c, \Sigma \H P^\infty)}\biggr) @>>>
\Sigma^2 C_1/A(\beta, Q_1) \to 0
\endmultline
$$
where
$$
H^*(\Sigma \C P^\infty_{-1}[1])/A(\Sigma y) \cong \F_p\{\Sigma y^k \mid
k \equiv 1 \bmod p-1, k \ge 1; k \ne p^e, e \ge 0\}
$$
and $C_1 \subset A$ is spanned over $\F_p$ by all admissible monomials
in $A$ except $1$ and the $P^I$ for $I = (p^e, p^{e-1}, \dots, p, 1)$
with $e\ge0$.
\endexample

\remark{Remark 5.6}
(a)\qua
The $A$-module $\Sigma^{-2} C/A(\beta, Q_1)$ can be shown to split off
from $H^*(Wh(*)/(\Sigma c, \Sigma \H P^\infty))$ by considering the
lower cofiber sequence in~3.6.

(b)\qua
For $p\ge5$ the extension of $\Sigma^2 C_1/A(\beta, Q_1)$ by
$H^*(\Sigma \C P^\infty_{-1}[1])/A(\Sigma y)$ is not split.  By~4.9 the
bottom $p$-torsion homotopy of $Wh(*)$ is $\Z/p$ in degree $4p-2$,
which implies that there is a nontrivial mod~$p$ Bockstein relating
the bottom classes $\Sigma^2 P^2$ and $\Sigma y^{2p-1}$ of these two
$A$-modules, respectively.
\endremark

\head Applications to automorphism spaces \endhead

We now recall the relation between Whitehead spectra, smooth concordance
spaces and diffeomorphism groups, to allow us to formulate a geometric
interpretation of our calculations.

\subhead Spaces of concordances and $h$-cobordisms \endsubhead
Let $M$ be a compact smooth $n$-manifold, possibly with corners,
and let $I = [0, 1]$ be the unit interval.  To study
the automorphism space $DIFF(M)$ of self-diffeomorphisms of $M$ relative
to the boundary $\partial M$, one is led to study the concordance space
$$
C(M) = DIFF(M \times I, M \times 1)
$$
of smooth concordances on $M$, also known as the pseudo-isotopy space
of $M$ \cite{Ig}.  This equals the space of self-diffeomorphisms $\psi$ of the
cylinder $M \times I$ relative to the part $\partial M \times I \cup M
\times 0$ of the boundary.  Both $DIFF(M)$ and $C(M)$ can be viewed as
topological or simplicial groups, and there is a fiber sequence
$$
DIFF(M \times I) @>>> C(M) @>r>> DIFF(M)
\tag 6.1
$$
where $r$ restricts a concordance $\psi$ to the upper end $M \times 1$
of the cylinder.

Let $J = [0, \infty)$.  The smooth $h$-cobordism space $H(M)$ of $M$
[\Cite{Wa3}, \S1] is the space of smooth codimension zero submanifolds $W
\subset M \times J$ that are $h$-cobordisms with $M = M \times 0$ at one
end, relative to the trivial $h$-cobordism $\partial M \times I$.  There
is a fibration over $H(M)$ with $C(M)$ as fiber and the contractible space
of collars on $M \times 0$ in $M \times J$ as total space.  Hence $H(M)$
is a non-connective delooping of $C(M)$, i.e., $C(M) \simeq \Omega H(M)$.
The homotopy types of the diffeomorphism group $DIFF(M)$, the concordance
space $C(M)$ and the $h$-cobordism space $H(M)$ are of intrinsic interest
in geometric topology.

There are stabilization maps $\sigma \: C(M) \to C(I \times M)$ and
$\sigma \: H(M) \to H(I \times M)$.  By Igusa's stability theorem
\cite{Ig}, the former map is at least $k$-connected when $n \ge
\max\{2k+7, 3k+4\}$.  Then this is also a lower bound for the connectivity
of the canonical map
$$
\Sigma \: C(M) @>>> \eusm C(M) = \hocolim_\ell C(I^\ell \times M)
$$
to the mapping telescope of the stabilization map $\sigma$ repeated
infinitely often.  We call $\eusm C(M)$ the stable concordance space
of $M$, and call the connectivity of $\Sigma \: C(M) \to \eusm C(M)$
the concordance stable range of $M$.  Likewise there is a stable
$h$-cobordism space $\eusm H(M) = \hocolim_\ell H(I^\ell \times M)$, and
$\eusm C(M) \simeq \Omega \eusm H(M)$.  The connectivity of the map $H(M)
\to \eusm H(M)$ is one more than the concordance stable range of $M$.

\subhead The stable parametrized $h$-cobordism theorem \endsubhead
Waldhausen proved in \cite{Wa6} that when $X = M$ is a compact smooth
manifold there is a homotopy equivalence
$$
\eusm H(M) \simeq \Omega \Omega^\infty Wh(M) \,,
\tag 6.2
$$
i.e., that the Whitehead space $\Omega^\infty Wh(M)$ of $M$ is a delooping
of the stable $h$-cobordism space $\eusm H(M)$ of $M$.  This stable
parametrized $h$-cobordism theorem is the fundamental result linking
algebraic $K$-theory of spaces to concordance theory.  At the level of
$\pi_0$ it recovers the (stable) $h$- and $s$-cobordism theorems of
Smale, Barden, Mazur and Stallings.  Waldhausen's theorem includes in
particular the assertion that the stable $h$-cobordism space $\eusm H(M)$
and the stable concordance space $\eusm C(M)$ are infinite loop spaces.

The functor $X \mapsto A(X)$ preserves connectivity of mappings, in
the sense that if $X \to Y$ is a $k$-connected map with $k\ge2$ then
$A(X) \to A(Y)$ is also $k$-connected [\Cite{Wa1}, 1.1], [\Cite{BM1}, 10.9].
It follows that $Wh(M)$, $\eusm H(M)$ and $\eusm C(M)$ take $k$-connected
maps to $k$-, $(k-1)$- and $(k-2)$-connected maps, respectively,
for $k\ge2$.

Let $\pi = \pi_1(M)$ be the fundamental group of $X = M$.  The classifying
map $M \to B\pi$ for the universal covering of $M$ is $k$-connected for
some $k\ge2$, so also $A(M) \to A(B\pi)$ is $k$-connected.  Let $R =
\Z[\pi]$.  Then the linearization map $L \: A(B\pi) \to K(R)$ is a rational
equivalence by [\Cite{Wa1}, 2.2].  Hence rational information about $K(R)$
gives rational information about $A(M)$ up to degree~$k$, and about
$\eusm C(M)$ up to degree~$k-2$, which in turn agrees with $C(M)$ in the
concordance stable range.

For example, Farrell and Hsiang \cite{FH} show that $\pi_m C(D^n)
\otimes \Q$ has rank~1 in all degrees $m \equiv 3 \bmod 4$, and rank~0
in other degrees, for $n$ sufficiently large with respect to $m$.
From this they deduce that $\pi_m DIFF(D^n) \otimes \Q$
has rank~1 for $m \equiv 3 \bmod 4$ and $n$ odd, and rank~0 otherwise,
always assuming that $m$ is in the concordance stable range for $D^n$.

For $\pi$ a finite group, $A(X)$ and $Wh(X)$ are of finite type by
theorems of Dwyer \cite{Dw} and Betley \cite{Be}, so the integral homotopy
type is determined by the rational homotopy type and the $p$-adic homotopy
type for all primes $p$.  Therefore our results on the $p$-adic homotopy
type of $Wh(*)$ have following application:

\proclaim{Theorem 6.3}
Assume~3.2.

{\rm(a)}\qua
Suppose $p\ge5$ and let $M$ be a $(4p-2)$-connected compact smooth manifold
whose concordance stable range exceeds $(4p-4)$, e.g., an $n$-manifold
with $n \ge 12p-5$.  Then the first $p$-torsion in the homotopy of
the smooth concordance space $C(M)$, and in the homotopy of the smooth
$h$-cobordism space $H(M)$, is
$$
\pi_{4p-4} C(M)_{(p)} \cong \pi_{4p-3} H(M)_{(p)} \cong \Z/p \,.
$$

{\rm(b)}\qua
Suppose $p=3$ and let $M$ be an $11$-connected compact smooth manifold
whose concordance stable range exceeds $9$, e.g., an $n$-manifold
with $n \ge 34$.  Then the first $3$-torsion in the homotopy of
the smooth concordance space $C(M)$, and in the homotopy of the smooth
$h$-cobordism space $H(M)$, is
$$
\pi_9 C(M)_{(3)} \cong \pi_{10} H(M)_{(3)} \cong \Z/3 \,.
$$
\endproclaim

\demo{Proof}
The first $p$-torsion in $\pi_* Wh(*)$ is $\Z/p$ in degree $* = 4p-2$ for
$p\ge5$, and $\Z/3\{\Sigma\beta_1\}$ in degree~$* = 11$ for $p=3$, and
$\pi_* Wh(*)$ is finite in all of these degrees.  When $M$
is $(4p-2)$-connected, resp{.} 11-connected, the map $\pi_* Wh(M) \to
\pi_* Wh(*)$ is an isomorphism in this degree.  And $\pi_{*-2} \eusm C(M)
\cong \pi_{*-1} \eusm H(M) \cong \pi_* Wh(M)$.  So if the concordance
stable range is at least $(4p-3)$, resp{.} 10, also $\pi_{*-2} C(M)
\cong \pi_{*-2} \eusm C(M)$ and $\pi_{*-1} H(M) \cong \pi_{*-1} \eusm
H(M)$ in this degree.
\qed
\enddemo

Similar statements may of course be given for when the subsequent torsion
groups in $\pi_* Wh(*)$ agree with $\pi_{*-2} C(M)$ and $\pi_{*-1} H(M)$,
under stronger connectivity and dimension hypotheses.

By [\Cite{KR}, 1.4] there is a summand $\Z/p$ in $\pi_{4p-4} C(M)$ for
any $p\ge5$, regular or not, but we need~3.2 to show that this is the
first $p$-torsion in $\pi_* C(M)$.

\proclaim{Theorem 6.4}
Assume~3.2.

{\rm(a)}\qua Suppose $p\ge5$ and let $M = D^n$ with $n \ge 12p-5$.
Then $\pi_{4p-4} DIFF(D^{n+1})$ or $\pi_{4p-4} DIFF(D^n)$ contains
an element of order~$p$.

{\rm(b)}\qua Suppose $p=3$ and let $M = D^n$ with $n \ge 34$.
Then $\pi_9 DIFF(D^{n+1})$ or $\pi_9 DIFF(D^n)$ contains an
element of order~$3$.
\endproclaim

\demo{Proof}
Consider the exact sequence in homotopy induced by~(6.1), with $D^n
\times I \cong D^{n+1}$.  A $\Z/p$ in $\pi_m C(D^n)$ either comes from
$\pi_m DIFF(D^{n+1})$, which known to be finite in these cases
by \cite{FH}, or maps to $\pi_m DIFF(D^n)$.
\qed
\enddemo


\newpage

\Refs

\ref \key {Ad} \by J\,F Adams
\paper Lectures on generalised cohomology
\inbook Category Theory, Hom\-ology Theory and their Applications
\jour Proc. Conf. Seattle Res. Center Battelle Mem. Inst. 1968,
\vol 3 \yr 1969 \pages 1--138
\endref

\ref \key {AP} \by J\,F Adams \by S\,B Priddy
\paper Uniqueness of $BSO$ 
\jour Math. Proc. Cambridge Philos. Soc.
\vol 80 \yr 1976 \pages 475--509
\endref

\ref \key {Be} \by S Betley
\paper On the homotopy groups of $A(X)$ 
\jour Proc. Amer. Math. Soc.
\vol 98 \yr 1986 \pages 495--498
\endref

\ref \key {BK} \by S Bloch \by K Kato
\paper $p$-adic {\'e}tale cohomology
\jour Publ. Math. Inst. Hautes {\'E}tud. Sci.
\vol 63 \yr 1986 \pages 107--152
\endref

\ref \key {BHM} \by M B{\"o}kstedt\by W-C Hsiang \by I Madsen
\paper The cyclotomic trace and algebraic $K$-theory of spaces
\jour Invent. Math. \vol 111 \yr 1993 \pages 865--940
\endref

\ref \key {BM1} \by M B{\"o}kstedt \by I Madsen
\paper Topological cyclic homology of the integers
\jour Asterisque \vol 226 \yr 1994 \pages 57--143
\endref

\ref \key {BM2} \by M B{\"o}kstedt \by I Madsen
\paper Algebraic $K$-theory of local number fields: the unramified case
\inbook Prospects in topology, Princeton, NJ, 1994
\bookinfo Ann. of Math. Studies \vol 138, \publ Princeton University Press
\yr 1995 \pages 28--57
\endref

\ref \key {Bo} \by A\,K Bousfield
\paper The localization of spectra with respect to homology
\jour Topology \vol 18 \yr 1979 \pages 257--281
\endref

\ref \key {Du} \by B\,I Dundas
\paper Relative $K$-theory and topological cyclic homology
\jour Acta Math. \vol 179 \yr 1997 \pages 223--242
\endref

\ref \key {Dw} \by W\,G Dwyer
\paper Twisted homological stability for general linear groups
\jour Ann. of Math. \vol 111 \yr 1980 \pages 239--251
\endref

\ref \key {DF1} \by W\,G Dwyer \by E\,M Friedlander
\paper Algebraic and {\'e}tale $K$-theory
\jour Trans. AMS \vol 292 \yr 1985 \pages 247--280
\endref

\ref \key {DF2} \by W\,G Dwyer \by E\,M Friedlander
\paper Topological models for arithmetic
\jour Topology \vol 33 \yr 1994 \pages 1--24
\endref

\ref \key {DM} \by W\,G Dwyer \by S\,A Mitchell
\paper On the $K$-theory spectrum of a ring of algebraic integers
\jour $K$-Theory \vol 14 \yr 1998 \pages 201--263
\endref

\ref \key {FH} \by F\,T Farrell \by W\,C Hsiang
\paper On the rational homotopy groups of the diffeomorphism groups of
discs, spheres and aspherical manifolds
\inbook Algebr. geom. Topol. Stanford/Calif. 1976
\bookinfo Proc. Symp. Pure Math. \vol 32, Part 1
\yr 1978 \pages 325--337
\endref

\ref \key {HM} \by L Hesselholt \by I Madsen
\paper On the $K$-theory of finite algebras over Witt vectors of perfect fields
\jour Topology \vol 36 \yr 1997 \pages 29--101
\endref

\ref \key {Ho} \by R\,T Hoobler
\paper When is $\Br(X)=\Br'(X)$?
\inbook Brauer groups in ring theory and algebraic geometry,
Proc. Antwerp 1981
\inbook Lecture Notes in Math. \vol 917 \publ Springer-Verlag
\yr 1982 \pages 231--244
\endref

\ref \key {Ig} \by K Igusa
\paper The stability theorem for smooth pseudoisotopies
\jour $K$-Theory \vol 2 \yr 1988 \pages 1--355 
\endref

\ref \key {KR} \by J.R Klein \by J Rognes
\paper The fiber of the linearization map $A(*) \to K(\Z)$
\jour Topology \vol 36 \yr 1997 \pages 829--848
\endref

\ref \key {Kn} \by K-H Knapp
\paper $Im(J)$-theory for torsion-free spaces
The complex projective space as an example
\paperinfo Revised version of Habilitationsschrift Bonn 1979,
in preparation
\endref

\ref \key {Li1} \by S Lichtenbaum
\paper On the values of zeta and $L$-functions,~I
\jour Annals of Math. \vol 96 \yr 1972 \pages 338--360
\endref

\ref
\key {Li2} \by S Lichtenbaum
\paper Values of zeta functions, {\'e}tale cohomology,
        and algebraic $K$-theory
\inbook Algebraic $K$-theory, II: ``Classical'' algebraic $K$-theory
and connections with arithmetic
\bookinfo Lecture Notes in Math. \vol 342 \publ Springer-Verlag
\yr 1973 \pages 489--501
\endref

\ref \key {MM} \by I Madsen \by R\,J Milgram
\book The classifying spaces for surgery and cobordism of manifolds
\bookinfo Annals of Mathematics Studies \vol 92,
\publ Princeton University Press
\yr 1979
\endref

\ref \key {MS} \by I Madsen \by C Schlichtkrull
\paper The circle transfer and $K$-theory
\eds Grove, Karsten et al
\inbook Geometry and topology, Aarhus
Proceedings of the conference on geometry and topology,
Aarhus, Denmark, August 10--16, 1998
\bookinfo Contemp. Math.  \vol 258 \publ American Mathematical Society
\yr 2000 \pages 307--328
\endref

\ref \key {MST} \by I Madsen\by V Snaith \by J Tornehave
\paper Infinite loop maps in geometric topology
\jour Math. Proc. Camb. Phil. Soc. \vol 81 \yr 1977
\pages 399--429 \endref

\ref \key {MT} \by I Madsen \by U Tillmann
\paper The stable mapping class group and $Q(\C P^\infty_+)$
\jour Invent. Math. \vol 145 \yr 2001 \pages 509--544
\endref

\ref \key{MQR} \by J\,P May\by F Quinn \by N Ray \by
        {\rm with contributions by} J Tornehave
\book $E_\infty$ ring spaces and $E_\infty$ ring spectra
\bookinfo Lecture Notes in Math. \vol 577 \publ Springer-Verlag
\yr 1977
\endref

\ref \key{Mil1} \by J\,S Milne
\book {\'E}tale cohomology
\bookinfo Princeton Mathematical Series \vol 33,
\publ Princeton University Press \yr 1980
\endref

\ref \key{Mil2} \by J\,S Milne
\book Arithmetic duality theorems 
\bookinfo Perspectives in Mathematics \vol 1
\publ Academic Press, Inc. \yr 1986
\endref

\ref \key{Mis} \by G Mislin
\paper Localization with respect to $K$-theory
\jour J. Pure Appl. Algebra \vol 10 \yr 1977/78 \pages 201--213
\endref

\ref \key{Mit} \by S\,A Mitchell
\paper On $p$-adic topological $K$-theory
\eds P\,G Goerss  et al
\inbook Algebraic $K$-theory and algebraic topology, Proceedings of the
NATO Advanced Study Institute, Lake Louise, Alberta, Canada, December
12--16, 1991
\bookinfo NATO ASI Ser. Ser. C, Math. Phys. Sci.
\vol 407 \publ Kluwer Academic Publishers
\yr 1993 \pages 197--204
\endref

\ref \key{Mo} \by R\,E Mosher
\paper Some stable homotopy of complex projective space
\jour Topology \vol 7 \yr 1968 \pages 179--193 
\endref

\ref \key{Mu} \by J Mukai
\paper The $S^1$-transfer map and homotopy groups of suspended complex
projective spaces
\jour Math. J. Okayama Univ. \vol 24 \yr 1982 \pages 179--200 
\endref

\ref \key{Qu1} \by D Quillen
\paper Higher algebraic $K$-theory. I
\inbook Algebraic $K$-Theory, I: Higher $K$-theories
\bookinfo Lecture Notes Math. \vol 341 \publ Springer-Verlag
\yr 1973 \pages 85--147 
\endref

\ref \key{Qu2} \by D Quillen
\paper Finite generation of the groups $K_i$ of rings of
algebraic integers
\inbook Algebraic $K$-Theory, I: Higher $K$-theories
\bookinfo Lecture Notes in Math. \vol 341 \publ Springer-Verlag
\yr 1973 \pages 179--198
\endref

\ref \key{Qu3} \by D Quillen
\paper Higher algebraic $K$-theory
\inbook Proc. Intern. Congress Math. Vancouver, 1974 \vol I
\publ Canad.~Math.~Soc. \yr 1975 \pages 171--176
\endref

\ref \key{Qu4} \by D Quillen
\paper Letter from Quillen to Milnor on
	$\operatorname{Im}(\pi_i O \to \pi_i^s \to K_i\Z)$
\inbook Algebraic $K$-theory, Proc. Conf. Northwestern Univ. Evanston,
Ill. 1976
\bookinfo Lecture Notes in Math. \vol 551 \publ Springer-Verlag
\yr 1976 \pages 182--188 
\endref

\ref \key{Ra} \by D\,C Ravenel
\book Complex cobordism and stable homotopy groups of spheres
\bookinfo Pure and Applied Math. \vol 121 \publ Academic Press \yr 1986
\endref

\ref \key{Ro} \by J Rognes
\paper Two-primary algebraic $K$-theory of pointed spaces
\jour Topology \vol 41 \yr 2002 \pages 873--926
\endref

\ref \key{RW} \by J Rognes \by C Weibel
\paper Two--primary algebraic $K$-theory of rings of integers in number fields
\jour J. Am. Math. Soc. \vol 13 \yr 2000 \pages 1--54 
\endref

\ref \key{St} \by N\,E Steenrod
\book Cohomology operations
\bookinfo Annals of Mathematics Studies \vol 50,
\publ Princeton University Press \yr 1962
\endref

\ref \key{SV} \by A Suslin \by V Voevodsky
\paper Bloch--Kato conjecture and motivic cohomology with finite coefficients
\eds B\,B Gordon et al
\inbook The arithmetic and geometry of algebraic cycles
\bookinfo NATO ASI Ser. Ser. C, Math. Phys. Sci.
\vol 548 \publ Kluwer Academic Publishers
\yr 2000 \pages 117--189
\endref

\ref \key{Ta} \by J Tate
\paper Duality theorems in Galois cohomology over number fields
\inbook Proc. Intern. Congress Math. Stockholm, 1962
\publ Inst. Mittag--Leffler \yr 1963 \pages 234--241
\endref

\ref \key{To} \by H Toda
\paper A topological proof of theorems of Bott and Borel--Hirzebruch for
homotopy groups of unitary groups
\jour Mem. Coll. Sci. Univ. Kyoto, Ser. A
\vol 32 \yr 1959 \pages 103--119
\endref

\ref \key{Vo1} \by V Voevodsky
\paper The Milnor Conjecture
\paperinfo preprint \yr 1996
\endref

\ref \key{Vo2} \by V Voevodsky
\paper On $2$-torsion in motivic cohomology
\paperinfo preprint \yr 2001
\endref

\ref \key{Wa1} \by F Waldhausen
\paper Algebraic $K$-theory of topological spaces. I
\inbook Algebr. geom. Topol. Stanford/Calif. 1976
\bookinfo Proc. Symp. Pure Math. \vol 32, Part 1
\yr 1978 \pages 35--60
\endref

\ref \key{Wa2} \by F Waldhausen
\paper Algebraic $K$-theory of topological spaces, II
\inbook Algebraic topology, Proc. Symp. Aarhus 1978
\bookinfo Lecture Notes in Math. \vol 763 \publ Springer-Verlag
\yr 1979 \pages 356--394
\endref

\ref \key{Wa3} \by F Waldhausen
\paper Algebraic $K$-theory of spaces, a manifold approach
\inbook Current trends in algebraic topology, Semin. London/Ont. 1981
\bookinfo CMS Conf. Proc. \vol 2, Part 1 \yr 1982 \pages 141--184
\endref

\ref \key{Wa4} \by F Waldhausen
\paper Algebraic $K$-theory of spaces
\inbook Algebraic and geometric topology,
Proc. Conf. New Brunswick/USA 1983
\bookinfo Lecture Notes in Math. \vol 1126 \publ Springer-Verlag
\yr 1985 \pages 318--419 
\endref

\ref \key{Wa5} \by F Waldhausen
\paper Algebraic $K$-theory of spaces, concordance, and stable homotopy theory
\inbook Algebraic topology and algebraic $K$-theory, Proc. Conf.
Princeton, NJ (USA)
\bookinfo Ann. Math. Stud. \vol 113 \yr 1987 \pages 392--417
\endref

\ref \key{Wa6} \by F Waldhausen {\rm et al}
\paper The stable parametrized $h$-cobordism theorem
\paperinfo in preparation
\endref

\endRefs
\enddocument